\documentclass[11pt,reqno]{amsart}

\usepackage{color}
\usepackage{collectbox}
\usepackage[all,cmtip,color]{xy}
\usepackage{amsthm}
\usepackage{amssymb}
\usepackage{latexsym}
\usepackage{multicol}
\usepackage{verbatim,enumerate}
\usepackage{hyperref}
\usepackage{nicematrix}
\usepackage{amsmath, amscd}
\usepackage{mathrsfs}
\usepackage[all,cmtip]{xy}
\usepackage{soul}

\advance\textwidth by 1.2in \advance\oddsidemargin by -.6in \advance\evensidemargin by -.6in
\usepackage{tikz}
\usetikzlibrary{decorations.markings}
\tikzstyle{vertex}=[ circle, fill, draw, inner sep=0pt, minimum size=4pt,]
\tikzstyle{edge}= [thick]

\newtheorem*{lem}{Lemma}
\newtheorem*{prop}{Proposition}

\newtheorem*{ex}{Example}

\theoremstyle{definition} \newtheorem*{defn}{Definition}
\theoremstyle{definition}
\newtheorem{thm}{Theorem}
\newtheorem*{thm*}{Theorem}

\newtheorem*{rem}{Remark}

\newenvironment{pf}{\proof}{\endproof}
\newcounter{cnt}
\newenvironment{enumerit}{\begin{list}{{\hfill\rm(\roman{cnt})\hfill}}{%
\settowidth{\labelwidth}{{\rm(iv)}}\leftmargin=\labelwidth%
\advance\leftmargin by \labelsep\rightmargin=0pt\usecounter{cnt}}}{\end{list}} \makeatletter
\def\mydggeometry{\makeatletter\dg@YGRID=1\dg@XGRID=20\unitlength=0.003pt\makeatother}
\makeatother \theoremstyle{remark}

\numberwithin{equation}{section}

 \DeclareMathOperator{\Ht}{ht}



\newcommand{\wt}{\operatorname{wt}}

\newcommand{\nc}{\newcommand}
\newcommand{\rnc}{\renewcommand}

\nc{\cal}{\mathcal} \nc{\goth}{\mathfrak} \rnc{\bold}{\mathbf}

\renewcommand{\Bbb}{\mathbb}
\nc\bomega{{\mbox{\boldmath $\omega$}}} \nc\bpsi{{\mbox{\boldmath $\Psi$}}}
 \nc\sing{{\rm sing}}
 \nc\balpha{{\mbox{\boldmath $\alpha$}}}
 \nc\bbeta{{\mbox{\boldmath $\beta$}}}
 \nc\bpi{{\mbox{\boldmath $\pi$}}}
  \nc\bpis{{\mbox{\boldmath \scriptsize$\pi$}}}
 \nc\bullets{{\mbox{\scriptsize $\bullet$}}}
 \nc\bvarpis{{\mbox{\boldmath \scriptsize$\varpi$}}}
  \nc\bvarpi{{\mbox{\boldmath $\varpi$}}}

\nc\bepsilon{{\mbox{\boldmath $\epsilon$}}}

  \nc\bomegas{{\mbox{\boldmath\scriptsize $\omega$}}}
  \nc\bepsilons{{\mbox{\boldmath \scriptsize$\epsilon$}}}
\nc{\spi}{{\rm sp}}
\nc\hlien{\hat{\lie n}^+}
  \nc\btaus{{\mbox{\boldmath \scriptsize$\tau$}}}\nc\bxi{{\mbox{\boldmath $\xi$}}}
\nc\bmu{{\mbox{\boldmath $\mu$}}} \nc\bcN{{\mbox{\boldmath $\cal{N}$}}} \nc\bcm{{\mbox{\boldmath $\cal{M}$}}} \nc\blambda{{\mbox{\boldmath
$\lambda$}}}%\nc\mathbb Nu{{\mbox{\boldmath $\nu$}}}
\nc\btau{{\mbox{\boldmath
$\tau$}}}

\newcommand{\lie}[1]{\mathfrak{#1}}

\makeatletter
\def\section{\def\@secnumfont{\mdseries}\@startsection{section}{1}%
  \z@{.7\linespacing\@plus\linespacing}{.5\linespacing}%
  {\normalfont\scshape\centering}}
\def\subsection{\def\@secnumfont{\bfseries}\@startsection{subsection}{2}%
  {\parindent}{.5\linespacing\@plus.7\linespacing}{-.5em}%
  {\normalfont\bfseries}}
\makeatother

 \nc{\Hom}{\operatorname{Hom}}
  \nc{\mode}{\operatorname{mod}}
\nc{\End}{\operatorname{End}} \nc{\wh}[1]{\widehat{#1}} \nc{\Ext}{\operatorname{Ext}}
 \nc{\ch}{\operatorname{ch}} \nc{\ev}{\operatorname{ev}}
\nc{\Ob}{\operatorname{Ob}} \nc{\soc}{\operatorname{soc}} \nc{\rad}{\operatorname{rad}} \nc{\head}{\operatorname{head}}

 \nc{\Cal}{\cal} \nc{\Xp}[1]{X^+(#1)} \nc{\Xm}[1]{X^-(#1)}
\nc{\on}{\operatorname} \nc{\Z}{{\bold Z}} \nc{\J}{{\cal J}} \nc{\C}{{\bold C}} \nc{\Q}{{\bold Q}}

\nc{\N}{{\Bbb N}} \nc\boa{\bold a} \nc\bob{\bold b} \nc\boc{\bold c} \nc\bod{\bold d} \nc\boe{\bold e} \nc\bof{\bold f} \nc\bog{\bold g}
\nc\boh{\bold h} \nc\boi{\bold i} \nc\boj{\bold j} \nc\bok{\bold k} \nc\bol{\bold l} \nc\bom{\bold m} \nc\bon{\bold n} \nc\boo{\bold o}
\nc\bop{\bold p} \nc\boq{\bold q} \nc\bor{\bold r} \nc\bos{\bold s} \nc\boT{\bold t} \nc\boF{\bold F} \nc\bou{\bold u} \nc\bov{\bold v}
\nc\bow{\bold w} \nc\boz{\bold z} \nc\boy{\bold y} \nc\ba{\bold A} \nc\bb{\bold B} \nc\bc{\mathbb C} \nc\bd{\bold D} \nc\be{\bold E} \nc\bg{\bold
G} \nc\bh{\bold H} \nc\bi{\bold I} \nc\bj{\bold J} \nc\bk{\bold K} \nc\bl{\bold L} \nc\bm{\bold M}  \nc\bo{\bold O} \nc\bp{\bold
P} \nc\bq{\bold Q} \nc\br{\bold R} \nc\bs{\bold S} \nc\bt{\bold T} \nc\bu{\bold U} \nc\bv{\bold V} \nc\bw{\bold W} \nc\bx{\bold
x} \nc\KR{\bold{KR}} \nc\rk{\bold{rk}} \nc\het{\text{ht }}
\nc\bz{\mathbb Z}
\nc\bn{\mathbb N}
\nc\us{\underline \bos}
\nc\uS{ \bs_{{\rm alt}}}
\nc\pr{\rm pr}

\nc\ts{Z(\bos)}

\nc\toa{\tilde a} \nc\tob{\tilde b} \nc\toc{\tilde c} \nc\tod{\tilde d} \nc\toe{\tilde e} \nc\tof{\tilde f} \nc\tog{\tilde g} \nc\toh{\tilde h}
\nc\toi{\tilde i} \nc\toj{\tilde j} \nc\tok{\tilde k} \nc\tol{\tilde l} \nc\tom{\tilde m} \nc\ton{\tilde n} \nc\too{\tilde o} \nc\toq{\tilde q}
\nc\tor{\tilde r} \nc\tos{\tilde s} \nc\toT{\tilde t} \nc\tou{\tilde u} \nc\tov{\tilde v} \nc\tow{\tilde w} \nc\toz{\tilde z} \nc\woi{w_{\omega_i}}
\nc\chara{\operatorname{Char}}

\begin{document}

\title[Monoidal categorification from alternating snakes]{Monoidal categorification from alternating snakes}
\author{Matheus Brito}
\address{Departamento de Matematica, UFPR, Curitiba - PR - Brazil, 81530-015}
\email{mbrito@ufpr.br}
\thanks{M.B. was partially supported CNPq grant 405793/2023-5}
\author{Vyjayanthi Chari}
\address{Department of Mathematics, University of California, Riverside, 900 University Ave., Riverside, CA 92521, USA}
\email{chari@math.ucr.edu}
\thanks{V.C. was partially supported by a travel grant from the Simons Foundation.}

\begin{abstract} In a recent paper,  the authors introduced the notion of an alternating snake and a corresponding family of finite dimensional modules for the quantum affine algebra associated to $A_n$. We prove that under some restrictions, an alternating snake defines a canonical monoidal category. We prove that this category has finitely many prime objects. As a consequence we prove that the Grothendieck ring is isomorphic to the Grothendieck ring of the category $\mathscr C_\xi$  for a suitable height function. In particular it follows that the special family of alternating snakes provides a monoidal categorification of a cluster algebra of type $A_N$ for a suitable value of $N$.
\end{abstract}

\maketitle
\section{Introduction}
The study of finite-dimensional representations of quantum affine algebras has been widely studied  for over three decades due to its rich structure and deep connections to various fields,  including, algebraic geometry, combinatorics, and integrable system in mathematical physics. \\\\
A recent and important  breakthrough in this area was made by Hernandez and Leclerc \cite{HL10,HL13a}, who proved that certain  tensor subcategories of finite-dimensional representations provide monoidal categorifications of cluster algebras. In particular, they constructed such subcategories of representations of quantum affine algebras of type $A_n$ by imposing compatibility conditions derived from a height function defined on the index set for the simple roots. They showed that in the case of the bipartite height function, the corresponding Grothendieck ring was the  monoidal categorification of a cluster algebra of type $A_n$ with $n$ frozen variables. This result was later proved in \cite{BC19a} for an arbitrary height function. \\\\
Around the same time, Mukhin and Young in \cite{MY12} introduced the snake modules for quantum affine $A_n$. This family generalizes the Kirillov--Reshetikhin modules and more generally the minimal affinizations and  has very nice properties. Except in small cases, they are not of the kind occurring in the  work of Hernandez and Leclerc. However the authors of the current paper noticed that there was an interesting pattern in these two families  modules.
This led them to introduce in \cite{BC25a} the notion of an alternating snake module for the quantum affine algebra of type $A_n$. These modules were proved to be  real and a necessary and sufficient  condition for them to be prime was also given. Moreover, under mild conditions, an explicit alternating formula for their character was also given. \\\\
The main goal of the present paper is to make a connection between alternating snake modules, tensor subcategories and monoidal categorification. More precisely, it is shown that, under suitable conditions, a prime  alternating snake $\bos$ determines a canonical monoidal subcategory $\mathscr C_n(\bos)$ of finite-dimensional representations of the quantum affine algebra $\bu_q(\widehat{\lie{sl}}_{n+1})$. The associated Grothendieck ring is isomorphic to a cluster algebra of type $A_N$, with $N$ determined explicitly by the  alternating snake. A central feature of the  approach is that these results are proved by giving an algorithm to describe tensor product factorizations of an arbitrary irreducible module in the category.
This is done by making   systematic use of $\ell$-highest weight theory and and  the Kashiawara--Kim--Oh--Park invariant $\lie d$ arising from the denominators of normalized $R$-matrices \cite{KKOP}.\\\\
These factorization results show that the Grothendieck ring is generated by classes of prime objects and that all irreducible classes factor uniquely in terms of these generators. This allows one to prove that the Grothendieck ring is isomorphic to the one coming from a particular class of height functions. The isomorphism maps irreducible (prime)  objects to irreducible (prime)  objects. Hence the results of \cite{BC19a} now imply that the monoidal subcategory defined by this family of alternating snakes is in fact a monoidal categorification of a cluster algebra.
\\\\
The paper is organized as follows.  Section \ref{monoidal} begins with  the background material on quantum affine algebras. The main result is stated after introducing the appropriate family of alternating snakes and  the category $\mathscr C_n(\bos)$ associated to a prime alternating snake $\bos$. Section \ref{Cismon} has the proof parts (i) and (ii) of Theorem \ref{mainthm}; namely that $\mathscr C_n(\bos)$ is a monoidal category and gives a  finite set of prime objects in $\mathscr C_n(\bos)$.  Section \ref{allprime}   develops the prime factorization theory for irreducible modules in $\mathscr C_n(\bos)$ and constitutes the technical core of the paper. These results are then used in Section \ref{isoGr} to study isomorphisms between Grothendieck rings arising from different alternating snakes and in Section \ref{moncat} to prove that the resulting category yields a  monoidal categorification of a cluster algebras of type $A$.

\section{Main Results}\label{monoidal}
We begin by recalling  some essential definitions and results on the representation theory of the quantum loop algebra $\widehat\bu_n$ associated to $\lie{sl}_{n+1}$. 
  As usual $\mathbb C$ (resp. $\mathbb C^\times$, $\mathbb Z$, $\mathbb Z_+$, $\mathbb N$) will denote the set of complex numbers (resp. non-zero complex numbers, integers, non-negative integers, positive integers).  Given $\ell\in \mathbb N$ we denote by $\Sigma_\ell$ the symmetric group on $\ell$ letters.  Assume throughout that $q$ is a non-zero complex number and not a root of unity.

\subsection{The algebra $\widehat\bu_n$ and the category $\mathscr F_n$}\label{basicdef}  For $n\in\mathbb N$, let  $\widehat\bu_n$ be  the quantum loop algebra associated to $\lie{sl}_{n+1}(\mathbb C)$;  we refer the reader to \cite{CP95} for precise definitions. For our purposes, it is enough to recall that $\widehat\bu_n$ is a Hopf algebra with 
an infinite set of generators: $x_{i,s}^\pm$, $\phi^\pm_{i,s}$,  $1\le i\le n$ and $s\in\mathbb Z$. The subalgebra $\widehat\bu_n^0$   generated by the elements $\phi^\pm_{i,s}$, $1\le i\le n$,  $s\in\mathbb Z$, is commutative.
\\\\
We shall be interested in  certain subcategories of the tensor category of finite-dimensional representations of $\widehat\bu_n$. Recall that an irreducible module is said to be prime if it is not isomorphic to a tensor product of non-trivial representations. A module is said to be real if its tensor square is irreducible.
It is well known (see \cite{CP91,CP95}) that the isomorphism  classes of irreducible finite-dimensional representations  of $\widehat\bu_n$  are indexed by the multiplicative monoid consisting  $n$-tuple of polynomials (Drinfeld polynomials) in a variable $u$ with constant term 1.   \\\\
Clearly any finite-dimensional module can be written as a direct sum of generalized eigenspaces for the action of $\widehat\bu_n^0$. The non-zero generalized eigenspaces are called $\ell$-weight spaces and the eigenvalues are called the $\ell$-weights of the module. It was shown in \cite{FR99} that the  $\ell$-weights of a finite dimensional module are contained in  the monoid  of $n$-tuples of rational functions in the indeterminate $u$. 

\subsection{An alternative formulation}\label{altform}
  Let $\mathbb I_n$ be the set of intervals $[i,j]$ with $i,j\in\mathbb Z$ and $0\le j-i\le n+1$. For $r\ge 1$ let $\mathbb I_n^r$ be the set of ordered $r$-tuples of elements of $\mathbb I_n$. Given elements $\bos_1\in\mathbb I_n^{r_1}$ and $\bos_2\in\mathbb I_n^{r_2}$ we let $\bos_1\vee\bos_2$ be the element of $\mathbb I_n^{r_1+r_2}$ obtained by concatenation.\\\\
Define $\cal I_n^+$ (resp. $\cal I_n$) to  be the  free abelian monoid (resp. group) with identity $\bold 1$ and generators
 $\bomega_{i,j}$ with $[i,j]\in\mathbb I_n$ where we understand that $\bomega_{i,i}=\bomega_{i,i+n+1}=\bold 1$ for all $i\in \mathbb Z$.
 Define map $\mathbb I_n^r\to\cal I_n^+$ given by  $ \bos=([i_1,j_1],\cdots, [i_r,j_r])\mapsto \bomega_{\bos}=\bomega_{i_1,j_1}\cdots\bomega_{i_r, j_r}$. Clearly given $\bomega\in\cal I_n^+$ there exists a unique positive integer $r$ and a unique (upto permutation) element $\bos=([i_1,j_1],\cdots, [i_r, j_r])\in\mathbb I_n^r$ with $0<j_s-i_s<n+1$ for all $1\le s\le r$ such that $\bomega=\bomega_{\bos}$ and we set $\Ht\bomega=r$.\\\\
 Identifying $\bomega_{i,j}$ with the Drinfeld polynomial $(1,\cdots, (1-q^{i+j}u), 1,\cdots, 1)$ (here the polynomial of degree one appears in the $(j-i)$-th position) we define the subcategory $\mathscr F_n$ to be the full subcategory of finite-dimensional representations whose simple objects are indexed by elements of $\cal I_n^+$. For any object $V$ of $\mathscr F_n$ let  $[V]$ denote the corresponding element of the Grothendieck group $ \cal K(\mathscr F_n)$ and for  $\bomega\in\cal I_n^+$ let $V(\bomega)$ be  the unique (up to isomorphism) irreducible $\widehat\bu_n$-module.
\\\\ 
 It was proved in \cite{HL10} that $\mathscr F_{n}$ is a rigid monoidal category and  hence $\cal K(\mathscr F_{n})$ admits a ring structure. The results of \cite{FR99} show that this     ring is polynomial in the generators $[V(\bomega_{i,j})]$, $[i,j]\in\mathbb I_n$. Moreover, it has a  basis given by the classes of the simple objects.  \\\\
 The results of \cite{FR99} show that the set $\wt_\ell V$ of $\ell$-weights of an object $V$  of $\mathscr F_n$ can be regarded as a  subset of $\cal I_n$ and that
$$\wt_\ell (V\otimes V')= (\wt_\ell V)(\wt_\ell  V'),\ \ (V\otimes V')_{\bomegas}= \bigoplus_{\bomegas_1\in\cal I_n} V_{\bomegas_1}\otimes V'_{\bomegas\bomegas_1^{-1}}. $$
Setting $\wt_\ell^+ V= \wt_\ell V\cap \cal I_n^+$ we note that if $V(\bomega)$ is a Jordan--H\"older constituent of $V$ then $\bomega\in\wt^+_\ell V$.

\subsection{Alternating snakes} The notion of an alternating snake and an alternating snake module was introduced in \cite{BC25a} and we  now  recall the definition for a particular family of alternating snakes.
\subsubsection{Sink Source type} Given $\bos=([i_1,j_1], \cdots, [i_r,j_r])\in\mathbb I_n^r$, set \begin{gather*}
a_{\max}=\max\{a_s: 1\le s\le r\},\ \ a_{\min}=\min\{a_s: 1\le s\le r\},\ \ a\in\{i,j\},\\
\bos(p-1,\ell)=([i_p,j_p],\cdots,[i_\ell,j_\ell]),\ \ 1\le p\le \ell\le r.\end{gather*}
\begin{defn}\label{altsnakedef}
 Say that an element $\bos=([i_1,j_1], \cdots, [i_r,j_r])\in\mathbb I_n^r$ is a stable alternating snake of sink source type if:
 \begin{itemize}
     \item 
for $1\le s< p\le r$  either $i_s\ne i_p$ or $j_s\ne j_p$,
\item  for all $1\le s\le r$ there exists $\epsilon_s\in\{0,1\}$ such that for $1\leq p\leq r-1$ we have $\epsilon_p+\epsilon_{p+1}=1$ and $ a_{p+1-\epsilon_p}<a_{p+\epsilon_p}$ for $a\in\{i,j\}$.
\item for all $1\le p\le s-2\le  r-2$ we have $i_p\le i_s<j_s\le j_p$.
\end{itemize}
    \end{defn}
 
\noindent  Let $\uS^{s}$ be the union of the set  of stable  alternating snakes of sink source type in $\mathbb I_n$ for all $n\ge 1$.  
 Clearly $\bos\in\uS^s$  if and only if $\bos(p-1,\ell)\in\uS^s$ for all $1\le p< \ell\le r$. 
 \subsubsection{Connected and Prime} Say that $\bos\in\uS^s$ is connected  if $$i_{p+1-\epsilon_p}<i_{p+\epsilon_p}\le j_{p+1-\epsilon_p}<j_{p+\epsilon_p} \ \ {\rm and} \ \ j_{p+\epsilon_p}-i_{p+1-\epsilon_p}\leq n+1,\ \ 1\le p\le r-1.$$
 Say that $\bos\in\uS^s$ is prime if it is connected and $$
a_s\neq a_{s+2} \ \ {\rm for }\ \ a\in\{i,j\}, \ \ 1\leq s\leq r-2.$$
Clearly $\bos$ is connected (resp. prime) if and only if $\bos(p-1,\ell)$ is connected (resp. prime) for all $1\le p\le \ell\le r$.
\subsubsection{Prime factor} \begin{defn} \label{pfdef} 
 We say that $\bos(0,p)$, $1\le p<r$, is a prime factor of $\bos\in\uS^s$ if  $\bos(0,p)$ is prime and either 
\begin{itemize}
\item
$([i_p, j_p],[i_{p+1}, j_{p+1}])$ is not connected,
 \item or $a_{p-1}=a_{p+1}$ for some $a\in\{i,j\}$ and if $b\in\{i,j\}$ with $b\ne a$ then
 $ b_{p+1-2\epsilon_p}<b_{p-1+2\epsilon_p}$,
 \item or $a_p=a_{p+2}$ for some $a\in \{i,j\}$ and if $b\in \{i,j\}$ with $b\neq a$ then $b_{p+2-2\epsilon_p}<b_{p+2\epsilon_p}$. 
 \end{itemize}\end{defn}
\subsection{The set $\mathbb I_n(\bos)$ and $\bop\bor(\bos)$}
Suppose that $\bos=([i_1,j_1], [i_2,j_2])$ is connected. Set
$$
\bop\bor(\bos)=\{[i_1,j_1], [i_2,j_2]\},\ \ \bof\bor(\bos)=\{([i_1,j_1], [i_2,j_2]),\  [i_1,j_2], \ [i_2,j_1]\}.$$
If  $\bos\in\uS^s\cap \mathbb I_n^r$ is prime with $r\ge 3$ set
\begin{gather*}\tilde{\mathbb I}_n(\bos)=\{[i_p, j_\ell]\in\mathbb I_n: 1\leq p\leq r, \  j_\ell\in\{j_{p-1-\epsilon_p}, j_{p-\epsilon_p}, j_{p+1-\epsilon_p}, j_{p+2-\epsilon_p}\}\},\\
\mathbb I_n(\bos)=\tilde{\mathbb I}_n(\bos)\setminus\{[i,i], [i, n+1+i]: i\in\mathbb Z\}.
\end{gather*}
For $\varepsilon,\varepsilon'\in\{0,1\}$ and for  $0\le p+1<\ell\le r$ set
$$\bos_{\varepsilon,\varepsilon'}(p,\ell+1)=([i_{p+\epsilon_p}, j_{p+1-\epsilon_p}])^\varepsilon\vee \bos(p+1,\ell)\vee ([i_{\ell+1+\epsilon_\ell}, j_{\ell+2-\epsilon_\ell}])^{\varepsilon'},$$ where we understand that,
\begin{itemize}
 \item the first term does not exist if $p=-1,0$ or if $\varepsilon=0$,
 \item the third term does not exist if $\ell=r-1,r$ or if $\varepsilon'=0$.
\end{itemize}
 Let $\bop\bor(\bos)$ be the union of $\mathbb I_n(\bos)$ with the set of elements of the form $\bos_{\varepsilon,\varepsilon'}(p,\ell+1)$ which satisfy,
\begin{itemize} \item   $\varepsilon=1$ only if $i_{p+\epsilon_p}\ne i_{p+3}$ and $j_{p+1-\epsilon_p}\ne j_{p+3}$,
 \item  $\varepsilon'=1$ only if  $i_{\ell-1}\ne i_{\ell+1+\epsilon_\ell}$ and $j_{\ell-1}\ne j_{\ell+2-\epsilon_\ell}$.
 \end{itemize}

Let $\bof\bor(\bos)$ be the union of the following  subsets of $\mathbb I_n\cup\mathbb I_n^2$:
\begin{gather*}
[i_{\min}, j_{\max}],\ \ [i_{\max}, j_{\min}],\\
    \{([i_1,j_1], [i_{2+\epsilon_1}, j_{3-\epsilon_1}]),\ \ ([i_s, j_s],[i_{s+1-2\epsilon_s}, j_{s-1+2\epsilon_s}]),\  \ ([i_r,j_r],[i_{r-2+\epsilon_r}, j_{r-1-\epsilon_r}]):\ \ 1<s<r\},\\
    \{([i_{s-\epsilon_s}, j_{s-1+\epsilon_s}], [i_{s+1+\epsilon_s}, j_{s+2-\epsilon_s}]): 2\le s\le r-2\ \ {\rm and}\  i_{s-1}\ne i_{s+2}, \ \ j_{s-1}\ne j_{s+2}\}. 
\end{gather*}

\subsection{The main result}
 Let $\cal I_n^+(\bos)$ be the submonoid with identity of $\cal I_n^+$ generated by the elements $\bomega_{i_s, j_p}$ with $[i_s, j_p]\in\mathbb I_n(\bos)$. Let $\mathscr C_n(\bos)$ be the full subcategory of $\mathscr F_n$ consisting of modules $V$ whose Jordan--Holder components are of the form $V(\bomega)$ with  $\bomega\in \cal I_n^+(\bos)$.
 The following is the main result of the paper.
 \begin{thm}\label{mainthm} Let $\bos\in\uS^s$ be prime.
 \begin{enumerit}
     \item[(i)] The category $\mathscr C_n(\bos)$ is monoidal.
     \item[(ii)] Let $\bos'\in\bop\bor(\bos)\sqcup\bof\bor(\bos)$. The module $V(\bomega_{\bos'})$ is prime and real.
     \item[(iii)] Given $\bold 1 \neq \bomega\in\cal I_n^+(\bos)$ there exist $k\in\mathbb N$ and elements $\bos_1,\cdots,\bos_k$ in $\bop\bor(\bos)\sqcup\bof\bor(\bos)$ such that
     $$V(\bomega)\cong V(\bomega_{\bos_1})\otimes\cdots\otimes V(\bomega_{\bos_k}).$$
     \item[(iv)] Assume that $\bos$ is such that $j_{\max}-i_{\min}=n+1$ and $j_{\min}=i_{\max}$. The Grothendieck ring $\cal K(\mathscr C_n(\bos))$ is a monoidal categorification of a cluster algebra of type $A_N$ where $N=r+\#\{2\le s\le r-3: i_{s-1}\ne i_{s+2}\ \ {\rm and}\ \ j_{s-1}\ne j_{s+2}\}.$
 \end{enumerit}
     
 \end{thm}
 The theorem is proved in the rest of the paper.
 \section{Proof of Theorem \ref{mainthm}{\rm (i),(ii)}}\label{Cismon}
\subsection{An enumeration} The following was proved in \cite{BC25a} in a more general setting. We include a proof in this particular case for the reader's convenience.
 \begin{lem}\label{enumerate} Suppose that $\bos=([i_1,j_1],\cdots,[i_r, j_r])\in\uS^s$ is prime; then for all appropriate $1\le s\le r$ we have
 \begin{gather*}
i_{s+1-\epsilon_s}<i_{s+2\epsilon_s}\le i_{s+3-2\epsilon_s}<i_{s+2+2\epsilon_s},\ \ \ j_{s+\epsilon_s}>j_{s+2-2\epsilon_s}>j_{s+1+2\epsilon_s}\ge j_{s+4-2\epsilon_s}.
 \end{gather*}
 In particular,
 \begin{gather*}
 i_s<i_{r-\epsilon_r}, \ \  s\ne r-\epsilon_r,\ \  i_{2-\epsilon_1}<i_p,\ \ p\ne 2-\epsilon_1,\\
 j_s< j_{1+\epsilon_1},\ \ s\ne 1+\epsilon_1,\ \ j_{r-1+\epsilon_r}< j_p,\ \ p\ne r-1+\epsilon_r.
 \end{gather*}
 \end{lem}
 \begin{pf}
     Suppose that $\epsilon_s=0$ in which case $\epsilon_{s+2}=0$ as well. The definition of  $\uS^s$ gives $i_{s+1}<i_s\le i_{s+3}<i_{s+2}$. Since $\bos$ is prime we also get $j_s>j_{s+2}>j_{s+1}\ge j_{s+4}$. The proof when $\epsilon_s=1$ is identical.
 \end{pf}
 \begin{rem}
     It is helpful to make the enumeration explicit. If  $\epsilon_1=0$ the lemma asserts that
     \begin{gather*}
         i_2<i_1\le i_4<i_3\le i_6<i_5<\cdots,\ \ j_1>j_3>j_2\ge j_5>j_4\ge j_7>j_6>\cdots,
     \end{gather*} while if $\epsilon=1$ we have
     \begin{gather*}
         i_1<i_3<i_2\le i_5<i_4\le i_7<\cdots,\ \ j_2>j_1\ge j_4>j_6>j_5\ge j_8>\cdots .
     \end{gather*}
 \end{rem}
 Note that Lemma \ref{enumerate} gives
\begin{equation}\label{in}{\rm if} \ j_{\max}-i_{\min}= n+1,\ \  {\rm and}\ \ j_{\min}=i_{\max}\ \ {\rm then}\ \ \mathbb I_n(\bos)=\tilde{\mathbb I}_n(\bos)\setminus\{[i_{\max}, j_{\min}],\,  [i_{\min}, j_{\max}]\}.\end{equation}
 \subsection{The set $\mathbb I_n(\bos)$}\label{closed}
 We need the following property
of $\mathbb I_n(\bos)$.
\begin{lem}\label{min}
    Suppose that $([i_p,j_\ell], [i_m,j_k])$ is a connected pair of elements in $\mathbb I_n(\bos)$ with $i_m<i_p\le j_k<j_\ell$. Then $[i_p, j_k]$ and $[i_m, j_\ell]$ are in $\tilde{\mathbb I}_n(\bos)$. 
\end{lem}
 \begin{pf} The lemma is proved by induction on $r$ with induction obviously beginning when $r=1,2$. We prove the inductive step when $\epsilon_1=0$; the proof in the other case is a trivial modification. We shall use Remark \ref{enumerate} freely without mention.\\\\
Assume that  $r>2$ so that the result holds for $\mathbb I_n(\bos(1,r))\subset\mathbb I_n(\bos)$. In particular, the lemma follows if we prove it in the case when  $j_1\in\{j_\ell, j_k\}$ or when  $i_1\in\{i_m, i_p\}$. If $j_1\in \{j_\ell,j_k\}$ then under the hypothesis of the lemma we have  $j_\ell=j_1$ and hence the  definition of $\mathbb I_n(\bos)$ gives $i_p\in\{i_1,i_2\}$. Since $i_m<i_p$ we are forced to have $i_p=i_1$ and $i_m=i_2$. It follows from the definition of $\tilde{\mathbb I}_n(\bos)$ that $j_k\in\{j_2, j_3\}$. Since 
$\{[i_1, j_3], [i_2,j_1], [i_1,j_2]\}\subset\tilde{\mathbb I}_n(\bos)$ by definition, the lemma follows in this case.\\\\
Assume from now on that $j_1\notin\{j_\ell, j_k\}$. Then, we may further assume that $p=1$ or $m=1$
(since otherwise the pair $([i_p,j_\ell], [i_m,j_k])$ are intervals in $\mathbb I_n(\bos(1,r))$ for which we know the result).
Hence together with the hypothesis of the Lemma, we get the following implications:
\begin{gather*}
    p=1, m=2\implies j_\ell\in\{j_2,j_3\},\ \ j_k\in\{j_2,j_3\}\setminus \{j_\ell\}\implies j_\ell=j_3,\ \ j_k=j_2,\\
    m=1\implies j_k\in\{j_2,j_3\}  \implies j_\ell=j_3,\ \  j_k=j_2\implies i_p\in\{i_3,i_4\}.
\end{gather*}
If $p=1$ the lemma follows since $[i_1, j_2]$ and $[i_2, j_3]$ are elements of $\tilde{\mathbb I}_n(\bos)$ while if $m=1$ the results follows since $[i_1,j_3]$, $[i_3, j_2]$ and $[i_4,j_2]$ are all elements of $\tilde{\mathbb I}_n(\bos)$.
\end{pf}
\begin{rem}
    In particular the preceding Lemma gives the following alternative characterization of $\tilde{\mathbb I}_n(\bos)$. Namely it is  the minimal subset  of $\mathbb I_n$ which contains $[i_p,j_p]$ for all $1\le p\le r$ and satisfies the following: if $[i,j]$ and $[i',j']$ is a connected pair of intervals in the set, then 
    $[i,j']$ and $[i',j]$ are also  elements of the set.
\end{rem}

\subsection{}\label{weyl}
Given $\bomega\in\cal I_n^+$ write $\bomega=\bomega_{p_1,\ell_1}\cdots\bomega_{p_k,\ell_k}\in\cal I_n^+$ so that one of the following hold: for $1\le m\le s\le k$ 
\begin{itemize}
    \item either the pair  $([p_m,\ell_m],[p_s, \ell_s])$ is  not connected,
    \item or $p_m+\ell_m\ge p_s+\ell_s$.
\end{itemize}
Set,
$$W(\bomega)= V(\bomega_{p_1,\ell_1})\otimes\cdots\otimes V(\bomega_{p_k,\ell_k}).$$
Then it is known (see \cite{Ch01,VV02}) that $W(\bomega)$ is an $\ell$--highest weight module; i.e., it is generated by a vector $w_\bomegas$ which is annihilated by the generators $x_{i,m}^+$ with $i\in I$, $m\in\mathbb Z$ and the elements of $\widehat\bu_q^0$ act as scalars on $w_\bomegas$ where the scalars are controlled by $\bomega$. Moreover, the module $V(\bomega)$ is the unique irreducible quotient of $W(\bomega)$.\\\\
In general we shall say that a module $V$ is an $\ell$--highest weight module with $\ell$--highest weight $\bomega$ if it is a quotient of $W(\bomega)$. It is also known that if $V(\bomega_1)\otimes V(\bomega_2)$ and $V(\bomega_2)\otimes V(\bomega_1)$ are both $\ell$--highest weight then $$V(\bomega_1\bomega_2)\cong V(\bomega_1)\otimes V(\bomega_2)\cong V(\bomega_2)\otimes V(\bomega_1).$$

\subsection{Proof of Theorem \ref{mainthm}(i)} The following result is immediate from  \cite[Theorem 1]{BC25b}.
\begin{prop}\label{lwt}
    Suppose that $\bos\in\uS^s$ and that $\bpi\in\cal I_n^+(\bos)$. 
    Then $\wt^+_\ell W(\bpi)\subset\cal I_n^+(\bos)$.\hfill\qedsymbol
\end{prop}
\noindent To prove Theorem \ref{mainthm}(i), suppose that $V(\bomega_1)$ and $V(\bomega_2)$ are objects of $\mathscr C_n(\bos)$; in particular $\bomega_1, \bomega_2\in\cal I_n^+(\bos)$.
 Taking $\bpi=\bomega_1\bomega_2$ in Proposition \ref{lwt} we see that $\wt_\ell^+ W(\bpi)$ is contained in $\cal I_n^+(\bos)$. Since $\wt_\ell^+(V(\bomega_1)\otimes V(\bomega_2))\subset \wt^+_\ell W(\bpi)$ it follows that if  $V(\bomega)$ occurs in the Jordan--Holder series of $V(\bomega_1)\otimes V(\bomega_2)$ then $\bomega\in\cal I_n^+(\bos)$. This proves that $\mathscr C_n(\bos)$ is a monoidal category.
\subsection{Proof of Theorem \ref{mainthm}(ii)}
\subsubsection{} The following was proved in \cite{BC25a}.
\begin{prop}\label{primefactor}
 Suppose that $\bos'\in\uS^s$. If $\bos'$ is prime then $V(\bomega_{\bos'})$ is prime and real. If $\bos'$ is not prime, then there exists $1\le p<r$ such that $\bos'(0,p)$ is prime and $$V(\bomega_{\bos'})\cong V(\bomega_{\bos'(0,p)})\otimes V(\bomega_{\bos'(p,r)}).$$\hfill\qedsymbol   
\end{prop}
\noindent In view of the proposition it suffices to prove that if $\bos'\in\bop\bor(\bos)\sqcup\bof\bor(\bos)$ then $\bos'$ is a prime element of $\uS^s$. If $\bos'\in\bof\bor(\bos)$ then it is straightforward to see that the pairs of intervals are connected and hence that $\bos'$ is prime.
\subsubsection{} Suppose that $\bos'=\bos_{\varepsilon,\varepsilon'}(p,\ell+1)\in \bop\bor(\bos)$. If $\varepsilon=0=\varepsilon'=0$ then $\bos'=\bos(p+1,\ell)$ and it is immediate from the fact that $\bos$ is prime element of $\uS^s$  that $\bos'$ is also a prime element of $\uS^s$. The following Lemma establishes the result in the remaining cases.

 \begin{lem}\label{prss} Suppose that $\bos\in\uS^s$ is prime and let $0\le p+1 < \ell\le r$.
 Then,
 \begin{gather*}
 a_{p}\neq a_{p+3}   \implies  \bos_{1,0}(p,\ell+1)\ {\rm is \ a \ prime\ element\ of  } \  \uS^s,\\
 a_{\ell-1}\neq a_{\ell+2}\implies  \bos_{0,1}(p,\ell+1)\ {\rm is \ a \ prime\ element\ of  } \  \uS^s,\\
a_{p}\neq a_{p+3}\ \ {\rm and}\ \ a_{\ell-1}\neq a_{\ell+2}\implies  \bos_{1,1}(p,\ell+1)\ {\rm is \ a \ prime\ element\ of  }  \ \uS^s. 
 \end{gather*}
 \end{lem}

 \begin{pf}
 If $a_p\neq a_{p+3}$ for $a\in \{i,j\}$, then Lemma \ref{enumerate} gives that
 \begin{gather*}
  i_{p+\epsilon_{p}}<i_{p+3}<j_{p+3}<j_{p+1-\epsilon_{p}},\ \
 i_{p+2\epsilon_{p}}<i_{p+2-\epsilon_p}\le j_{p+1+\epsilon_p}<j_{p+2-2\epsilon_p},\\
 s>p+3 \implies i_{p+\epsilon_{p}}\leq i_{s}<j_{s}\leq j_{p+1-\epsilon_{p}}.
 \end{gather*} 
The preceding inequalities show that  $\bos_{1,0}(p,\ell+1)=[i_{p+\epsilon_{p}}, j_{p+1-\epsilon_{p}}]\vee \bos(p+1,\ell)$   is  a prime element of $\uS^s$. Similarly, if $a_{\ell-1}\neq a_{\ell+2}$, then Lemma \ref{enumerate} also  gives, \begin{gather*} i_{\ell-1}<i_{\ell+1+\epsilon_\ell}\le j_{\ell+2-\epsilon_\ell}<j_{\ell-1}, \ \ i_{\ell+1-\epsilon_\ell}< i_{\ell+2\epsilon_\ell}\le j_{\ell+2-2\epsilon_\ell}<j_{\ell+\epsilon_\ell},\\
 s<\ell-1 \implies i_s\leq i_{\ell+1+\epsilon_\ell}\le j_{\ell+2-\epsilon_\ell}\le j_s,
 \end{gather*}
and hence $\bos_{0,1}(p,\ell+1)=\bos(p+1,\ell)\vee([i_{\ell+1+\epsilon_\ell}, j_{\ell+2-\epsilon_\ell}])$ is a prime element of $\uS^s $. 
The proof of the third assertion is now immediate.
\end{pf}
\section{Proof of Theorem \ref{mainthm}{\rm (iii)}}\label{allprime}
We begin this section by recalling some results of \cite{KKOP} which will be needed for our work.

\subsection{} \label{kkop}In \cite{KKOP}, the authors defined  for   $\bomega_1, \bomega_2\in\cal I_{n}^+$ a  non-negative integer $\lie d(V(\bomega_1), V(\bomega_2))$ depending on $n$.  This integer corresponds to the order of the zero of $d_{\bomegas_1,\bomegas_2}(z)d_{\bomegas_2,\bomegas_1}(z)$, where $d_{\bomegas_1,\bomegas_2}(z)$ denotes the denominator of the normalized $R$-matrix of $V(\bomega_1)$ and $V(\bomega_2)$. 
We summarize certain important properties of $\lie d$ in the following proposition. Part (i) follows from the definition of $\lie d$, part (ii) is Corollary 3.17 of \cite {KKOP} and part (iii) is Proposition 4.2 of \cite{KKOP2}. Part  (iv) follows from Proposition 2.17 of \cite{KKOP22}.
Part (v) was proved in \cite{Naoi24}.

\begin{prop} Let $\bomega_1, \bomega_2\in\cal I^+_{n}$ and assume that  $V(\bomega_1)$ is a real $\widehat\bu_n$--module. Then,
\begin{enumerit}
    \item [(i)] $\lie d(V(\bomega_1), V(\bomega_2))=\lie d(V(\bomega_2), V(\bomega_1))$.
    \item[(ii)] $\lie d(V(\bomega_1), V(\bomega_2))=0$ if and only if $V(\bomega_1)\otimes V(\bomega_2)$ is irreducible.
    \item[(iii)] For all $\bomega_3\in\cal I_n^+$ we have $$\lie d(V(\bomega_1), V(\bomega_2\bomega_3))\le \lie d (V(\bomega_1), V(\bomega_2))+\lie d(V(\bomega_1), V(\bomega_3)).$$
    \item[(iv)] If $\lie d (V(\bomega_1), V(\bomega_2))\leq 1$ then $\lie d (V(\bomega_1\bomega_2), V(\bomega_1))=0$.
     \item[(v)] Suppose that $\bomega=\bomega_{i_1,j_1}\cdots\bomega_{i_k,j_k}$ is such that $i_s<i_{s+1}$ and $j_s<j_{s+1}$ for all $1\le s\le k-1$. Then $V(\bomega')\otimes V(\bomega)$ is irreducible if $\bomega(\bomega')^{-1}\in\cal I_n^+$.

    \end{enumerit}
    \end{prop}

\subsection{}
Part (iii) of Theorem \ref{mainthm} is immediate from the next proposition.

\begin{prop}\label{pd} Given $\bomega\in\cal I_n^+(\bos)$ there exists $\bos_1\in\bop\bor(\bos)\sqcup \bof\bor(\bos)$ such that $$V(\bomega)\cong  V(\bomega_{\bos_1})\otimes V(\bomega\bomega_{\bos_1}^{-1}).$$
\end{prop}
 The proof of the proposition proceeds by induction on $r$ and for a fixed $r$ by a further induction on $\Ht\bomega$ (see Section \ref{altform} for the definition). If $r=1$ there is nothing to prove since in that case we have $\bomega=\bomega_{i_1,j_1}^s$ for some $s\ge 1$ and $V(\bomega_{i_1,j_1})^{\otimes s}\cong V(\bomega_{i_1,j_1}^s)$. The proof of the inductive step is given in the rest of the section.
\\
\subsection{}\label{comments} We make some preliminary comments.
We shall use the enumeration of $\bos$ given in Lemma \ref{enumerate} repeatedly without mention. We shall also use some parts of Proposition \ref{kkop} very frequently also without mention. Thus, instead of proving that $V(\bomega_1)\otimes V(\bomega_2)$ is irreducible we will establish the equivalent formulation given in Proposition \ref{kkop}(ii).  Similarly, if we want to prove that $\lie d (V(\bomega_1\bomega_2), V( \bomega_3))=0$ then we shall just prove  that $\lie d (V(\bomega_s), V(\bomega_3))=0$ for $s=1,2$. \\\\
It is  straightforward to observe that if  $$ \mathbb I_n(\bos_{1,0}(1,r+1))=\{[i_{1+\epsilon_1},j_{2-\epsilon_1}], [i_{2+\epsilon_2}, j_{3-\epsilon_2}], [i_{2+2\epsilon_2}, j_{4-2\epsilon_2}] , [i_{1+2\epsilon_1},j_{3-2\epsilon_1}]\}\sqcup\mathbb I_n(\bos(2,r)).$$ Note that for all $\bos'\in\bop\bor(\bos)\sqcup\bof\bor(\bos)$ we have $\mathbb I_n(\bos')\subset\mathbb I_n(\bos)$ and so  $\cal I_n^+(\bos')$ is a submonoid of $\cal I_n^+(\bos)$. Hence, for the inductive step  we can assume that \begin{equation*}\label{indass} \bomega\notin\cal I_n^+(\bos(1,r)).\end{equation*} In particular we can and do write
\begin{equation}\label{omega}\bomega=\bomega_{i_1,j_1}^{a_1}\bomega_{i_{1+\epsilon_1}, j_{2-\epsilon_1}}^{a_2}\bomega_{i_{1+2\epsilon_1}, j_{3-2\epsilon_1}}^{a_3}\bomega_{i_{2-\epsilon_1}, j_{1+\epsilon_1}}^b\bomega',\ \ \ \ \bomega'\in\cal I_n^+(\bos(1,r))\end{equation}
where:
\begin{itemize}\item where $a_2=a_3=0$ if $\epsilon_1=0$ and  $i_1=i_4$ or $\epsilon_1=1$ and $j_1=j_4$,
\item $a_1+a_2+a_3+b>0$.\end{itemize}

Finally, we shall state all results for  $\epsilon_1\in\{0,1\}$, but for ease of notation only prove them when $\epsilon_1=0$. The proof is identical in the other case.

\subsection{} \label{ba30} Recall that  $i_{2-\epsilon_1}=i_{\min}$ and $j_{1+\epsilon_1}= j_{\max}$ and so $[i_{2-\epsilon_1}, j_{1+\epsilon_1}]$ is not connected to any  $[i_s,j_p]$ for all $1\le s,p\le  r$. Hence if $b\ge 1$  we can take $\bos_1=[i_{2-\epsilon_1}, j_{1+\epsilon_1}]$. 
  Next we deal with the case when $a_3\ge 1$ for which we need the following proposition.
\begin{prop}\label{red1}
    \begin{enumerit}
        \item[(i)] The interval $[i_{1+2\epsilon_1}, j_{3-2\epsilon_1}]$ and $[i_s, j_p]\in\mathbb I_n(\bos)$ are connected if and only if $i_s=i_2$ and $j_p=j_2$. 
        \item[(ii)] If $[i_2,j_2]$ and $[i_s,j_p]\in\mathbb I_n(
        \bos)$ are connected then $$[i_s,j_p]\in\{[i_1,j_1], [i_{1+2\epsilon_1}, j_{3-2\epsilon_1}], [i_3,j_3],[i_{4-\epsilon_1},j_{3+\epsilon_1}]\}.$$
     \item[(iii)]   
   We have      
$$\lie d (V(\bomega_{i_{1+2\epsilon_1},\, j_{3-2\epsilon_1}}\bomega_{i_2,j_2}), V(\bomega_{i_s,j_p}))=0,\ \ [i_s, j_p]\in\mathbb I_n(\bos).$$
   
    \end{enumerit}
\end{prop}
\begin{pf} Suppose that $\epsilon_1=0$. If $s
\ge 3$, then  by  the definition of $\mathbb I_n(\bos)$ we have $p\ge 2$. Hence $i_1<i_s<j_p\le j_3$ and so the intervals $[i_1,j_3]$ and $[i_s,j_p]$ are not connected. 
If  $s=2$ then by definition of $\mathbb I_n(\bos)$ we have $j_p\in\{j_1,j_2,j_3\}$. It follows that $[i_2,j_p]$ and $[i_1,j_3]$ are connected if and only if $j_p=j_2$ which proves (i).
\\\\
For part (ii) we must have $i_2<i_s\le j_2<j_p$ and hence we have  $j_p\in \{j_1,j_3\}$. If $j_p=j_1$ the definition of $\mathbb I_n(\bos)$ forces $i_s=i_1$; if $j_p=j_3$ the definition forces  $i_s\in\{i_1, i_3, i_4\}$ and part (ii) is proved.
\\\\
Using parts (i) and  (ii)  of the proposition we see that it suffices to prove (iii) when,
$$[i_s,j_p]\in\{[i_1,j_3], [i_2,j_2], [i_1,j_1], [i_3,j_3], [i_4,j_3]\}.$$ In the first two cases the assertion follows from Proposition \ref{kkop}(v). In the remaining three cases, notice that (after a possible reordering of the first and third term) we have that $([i_1,j_3], [i_2,j_2], [i_s,j_p])\in\uS^s$ and is not prime. It is easy to check using Definition \ref{pfdef} that $([i_1,j_3], [i_2,j_2])$ is a prime factor and  hence the result follows from Proposition \ref{primefactor}.\end{pf}
\noindent  Suppose that $\epsilon_1=0$ and $a_3\ge 1$ in \eqref{omega}.  Proposition \ref{red1}(i) gives $$\lie d(V(\bomega_{i_1,j_3}), V(\bomega\bomega_{i_1,j_3}^{-1}))=0\ \ \ {\rm if}\ \ \bomega'\bomega_{i_2,j_2}^{-1}\notin\cal I_n^+,$$ and so we can take $\bos_1=[i_1,j_3]$ if $\bomega'\bomega_{i_2,j_2}^{-1}\notin\cal I_n^+$. Otherwise, if  $\bomega'\bomega_{i_2,j_2}^{-1}\in\cal I_n^+$ then part (iii) of the proposition shows that we can take $\bos_1=([i_1,j_3], [i_2,j_2])$.
\subsection{} We now prove that the proposition holds if  $a_3=0=b$ and $\bomega'\bomega_{i_{3-\epsilon_2}j_{2+\epsilon_2}}^{-1}\in\cal I_n^+$. This relies on the following proposition.
\begin{prop} \label{i2j3}
\begin{enumerit}
    \item[(i)] Suppose that $[i_1,j_1]$ and $[i_s, j_p]\in \mathbb I_n(\bos)$ are connected. If $\epsilon_1=0$ (resp. $\epsilon_1=1$) then $i_s=i_2$ (resp. $j_p=j_2$) and $j_p\in\{j_2, j_3\}$ (resp. $i_s\in\{i_2,i_3\}$).
    \item[(ii)] The intervals $[i_{3-\epsilon_2},j_{2+\epsilon_2}]$ and $[i_s, j_p]$ are connected iff $i_s=i_1$ and $j_p=j_1$. In particular for $[i,j]\in\mathbb I_n(\bos)$ we have \begin{gather*}
    \lie d( V(\bomega_{i_{3-\epsilon_2},j_{i+\epsilon_2}}), V(\bomega_{i,j}))=0\ \ {\rm if} \ \ i\ne i_1\ \ {\rm or}\ \ j\ne j_1,\\  \lie d(V(\bomega_{i_1,j_1}\bomega_{i_{3-\epsilon_2},j_{2+\epsilon_2}}), V(\bomega_{i,j}))=0.\end{gather*}
\end{enumerit}
\end{prop}
\begin{pf} Suppose that $\epsilon_1=0$. By definition, $[i_2,j_p]\in\mathbb I_n(\bos)$ only if $j_p=j_1, j_2,j_3$.
For (i) we note that $j_p<j_1$  and hence we must have $i_s<i_1$ which implies $i_s=i_2$ and so $j_p=j_2,j_3$.\\\\  For part (ii) we note that since $i_2<i_s$ we must have $j_3<j_p$. This forces $j_p=j_1$ and since $[i_s,j_1]\in\mathbb I_n(\bos)$ only if $i_s=i_1, i_2$ we see that $i_s=i_1$.\\\\
The remaining statements in part (ii) are now immediate except for the second equality and the case $[i,j]=[i_2,j_2]$. Here we observe that  $([i_2,j_3],[i_1,j_1], [i_2,j_2])\in\uS^s$ is not prime and hence by Proposition \ref{primefactor}  we have that $$\lie d (V(\bomega_{i_1,j_1}\bomega_{i_2,j_3}), V(\bomega_{i_2,j_2}))=0.$$
\end{pf}

As a consequence of Proposition \ref{i2j3} and Proposition \ref{kkop}(iii),  we see that 
if $\bomega\bomega_{i_{3-\epsilon_2}, j_{2+\epsilon_2}}^{-1}\in\cal I_n^+$
then either
 $\bos_1=([i_1,j_1], [i_{3-\epsilon_2},j_{2+\epsilon_2}])$ or $\bos_1=[i_{3-\epsilon_2},j_{2+\epsilon_2}]$.
\\\\ 
Hence we are reduced to prove the proposition in the following case:
\begin{equation}\label{omega1}
\bomega=\bomega_{i_1,j_1}^{a_1}\bomega_{i_{1+\epsilon_1},j_{2-\epsilon_1}}^{a_2}\bomega_{i_2, j_2}^c\bomega',\ \ a_1+a_2>0,\end{equation} where 
\begin{equation}\label{rest1}\bomega'\bomega_{i,j}^{-1}\notin\cal I_n^+\ \ {\rm if}\ \ [i,j]\in\{[i_2,j_2], [i_{2-\epsilon_1},j_{1+\epsilon_1}], [i_{1+2\epsilon_1},j_{3-2\epsilon_1}],[i_{3-\epsilon_2},j_{2+\epsilon_2}]\}.\end{equation}
In particular, we have that if $[i_1,j_1]$ and $[i,j]$ are connected with $\bomega\bomega_{i,j}^{-1}\in\cal I_n^+$ then $[i,j]=[i_2, j_2]$.
If $a_1>c$ then writing $\bomega=\bomega_{i_1,j_1}^{a_1-c}(\bomega_{i_1,j_1}\bomega_{i_2, j_2})^c\bomega''$ we get 
$$\lie d (V(\bomega_{i_1,j_1}^{a_1-c}), V(\bomega_{i_1,j_1}^c\bomega_{i_2,j_2}^c\bomega''))\leq \lie d (V(\bomega_{i_1,j_1}^{a_1-c}), V(\bomega_{i_1,j_1}^c\bomega_{i_2,j_2}^c)) + \lie d (V(\bomega_{i_1,j_1}^{a_1-c}), V(\bomega''))=0,$$
where the first inequality follows from Proposition \ref{kkop}(iii) and the last equality follows from Proposition \ref{kkop}(v) and Proposition \ref{i2j3}(i); then we can take $\bos_1=[i_1,j_1]$.
\subsection{} \label{i1j2} We assume from now on that we are in the situation of \eqref{omega1}--\eqref{rest1} and that $a_1\le c$.
We   prove that the proposition holds if $a_2>0$. As usual we shall prove this when $\epsilon_1=0$ and so by the conditions following \eqref{omega}  we can and do assume that $i_1\ne i_4$. In particular, $\bos_{1,0}(1,r+1)\in\uS^s$ is prime by Lemma \ref{prss}. Since $\bomega_{i_1,j_2}^{a_2}\bomega'\in\cal I_n^+(\bos_{1,0}(1,r+1))$   the inductive hypothesis applies to this element.  Specifically, the hypothesis gives that there exist
\begin{itemize}\item integers $0\le d, d'\le a_2$ 
\item for $1\le k\le a_2-d-d'$, elements 
$\bos_{1,\varepsilon_k}(1,s_k+1)\in\bop\bor(\bos_{1,0}(1,r+1)),\ \ \varepsilon_k\in\{0,1\}$\end{itemize}
such that setting $$\bomega_3=\bomega_{i_1,j_2}^{a_2}\bomega'\left(\bomega_{\bos_{1,\varepsilon_1}(1,s_1+1)}\cdots \bomega_{\bos_{1,\varepsilon_k}(1,s_k+1)}\right)^{-1}\bomega_{i_4,j_3}^{-d'} ,
$$ we have 
\begin{equation}\label{omega3}V(\bomega_{i_1,j_2}^{a_2}\bomega')\cong V(\bomega_{i_1,j_2})^{\otimes d}\otimes V(\bomega_{i_1,j_2}\bomega_{i_4, j_3})^{\otimes d'}\bigotimes_{k=1}^{a_2-d-d'} V(\bomega_{\bos_{1,\varepsilon_k}(1,s_k+1)})\otimes V(\bomega_3).\end{equation}
\\\\
If $c=0$ then $a_1=0$ and so there is nothing to prove. The following proposition
is needed to prove the case when $c\ge a_1>0$.

\begin{prop} 
    For $p\in\{1,2\}$ we have,\begin{gather*}
    \lie d(V(\bomega_{i_p, j_p}), V(\bomega_{i_1,j_2}))=0=\lie d( V(\bomega_{i_p,j_p}), V(\bomega_{i_1,j_2}\bomega_{i_4,j_3})),\\ 
    \lie d (V(\bomega_{i_p, j_p}), V(\bomega_{\bos_{1,\varepsilon}(1,s)}))=0,\  \ \varepsilon\in \{0,1\}, \ \ \ 2\le s\le r+1.\end{gather*}\end{prop}
    \begin{pf} The first equality is trivial since for $p=1,2$ the intervals $[i_p,j_p]$ and $[i_1,j_2]$ are not connected. 
    The second holds when $p=1$ since $i_1\le i_4\le j_3<j_1$. If $p=2$ then note that $([i_2,j_2], [i_4,j_3],[i_1,j_2])\in\uS^s$ and is not prime. An application of Proposition \ref{primefactor} gives the result.
\\\\
The third equality follows from  Proposition \ref{i2j3}(i) if $p=1$. If $p=2$ then an application of Proposition \ref{red1}(ii) and Proposition \ref{kkop}(iii) shows that it is enough to prove that if $\epsilon_1=0$,
    $$\lie d(V(\bomega_{i_2,j_2}), V(\bomega_{i_1,j_2}\bomega_{i_3,j_3}))=0.$$ This again follows from Proposition \ref{primefactor} since $([i_2,j_2],[i_3,j_3],[i_1,j_2])\in\uS^s$ is not prime.\\
    \end{pf}
  Together with    an application of Proposition \ref{kkop}  the preceding proposition now gives 
    $$\lie d (V , V(\bomega_{i_1,j_1}^{a_1}\bomega_{i_2,j_2}^c\bomega_3))=0$$ where $V$ is one of the tensor factors appearing in \eqref{omega3} and hence we have
    \begin{equation}\label{a2>0case}
        V(\bomega)\cong V(\bomega_{i_1,j_2})^{\otimes d}\otimes V(\bomega_{i_1,j_2}\bomega_{i_4, j_3})^{\otimes d'}\bigotimes_{k=1}^{a_2-d-d'} V(\bomega_{\bos_{1,\varepsilon_k}(1,s_k+1)})\otimes V(\bomega_{i_1,j_1}^{a_1}\bomega_{i_2,j_2}^c\bomega_3).
    \end{equation}
    In particular Proposition \ref{pd} follows if $a_2>0$.
   
     \subsection{}\label{finalcase}To complete the proof of Proposition \ref{pd} in the case when $\epsilon_1=0$ we now have 
$$\bomega=\bomega_{i_1,j_1}^{a_1}\bomega_{i_2,j_2}^c\bomega', \ \ a\le c$$ where $$\bomega'\bomega_{i,j}^{-1}\notin\cal I_n^+,\ \ [i,j]\in\{[i_1, j_1], [i_1,j_2], [i_2,j_2], [i_2,j_1], [i_1,j_3], [i_2,j_3]\}.$$ 
   The inductive hypothesis on $r$ applies to $\bomega_{i_2,j_2}^c\bomega'$  and hence there exists 
   \begin{itemize}\item $(\ell_s,\varepsilon_s)$ with $\ell_s\ge 2$, $\varepsilon_s\in\{0,1\}$ for $1\le s\le c$,
   \item $\bomega''\in\cal I_n^+(\bos(1,r))$,
   \end{itemize} such that
\begin{gather}\label{indhyp2r}V(\bomega_{i_2,j_2}^c\bomega')\cong V(\bomega_{\bos_{0,\varepsilon_1}(0,\ell_1+1)})\otimes\cdots\otimes V(\bomega_{\bos_{0,\varepsilon_c}(0,\ell_c+1)})\otimes V(\bomega''),\ \ \ \ \bomega''\bomega_{i_2,j_2}^{-1}\notin\cal I_n^+.
\end{gather}
In particular for $ 1\le s,p\le c$, we have
\begin{equation}\label{omega''}\lie d(V(\bomega_{\bos_{0,\varepsilon_s}(0,\ell_s+1)}), V(\bomega_{\bos_{0,\varepsilon_p}(0,\ell_p+1)}))=0=\lie d(V(\bomega_{\bos_{0,\varepsilon_s}(0,\ell_s+1)}), V(\bomega'')).\end{equation}
Assume  without loss of generality 
 that $\ell_1,\cdots, \ell_c$ is ordered so that:
\begin{itemize}
\item if $\ell_s+\varepsilon_s$ is odd then  $\ell_{s+1}+\varepsilon_{s+1}$ is also odd and $2\ell_{s}+\varepsilon_{s}\ge  2\ell_{s+1}+\varepsilon_{s+1}$
 \item if $\ell_{s}+\varepsilon_{s}$ and $\ell_{s+1}+\varepsilon_{s+1}$ are both even then $2\ell_{s}+\varepsilon_{s}\leq 2\ell_{s+1}+\varepsilon_{s+1}$.
 \end{itemize}
The  final case of Proposition \ref{pd} follows if we prove that
$$V(\bomega)\cong \bigotimes_{s=1}^{a_1}V(\bomega_{\bos_{0,\varepsilon_s}(-1,\ell_s+1)})\bigotimes_{s=a_1+1}^cV(\bomega_{\bos_{0,\varepsilon_s}(0,\ell_s+1)})\otimes V(\bomega'').$$ 
\noindent For this, it suffices to prove that for $\, 1\le s\le s'\le a$ and $ a+1\le s''\le c$ we have
\begin{gather*}
    \lie d(V(\bomega''),V(\bomega_{\bos_{0,\varepsilon_s}(-1,\ell_s+1)}))=0=
     \lie d(V(\bomega_{\bos_{0,\varepsilon_{s}}(-1,\ell_{s}+1)}),V(\bomega_{\bos_{0,\varepsilon_{s''}}(0,\ell_{s''}+1)})),\\
     \lie d(V(\bomega_{\bos_{0,\varepsilon_{s}}(-1,\ell_{s}+1)}),  V(\bomega_{\bos_{0,\varepsilon_{s'}}(-1,\ell_{s'}+1)}))=0.
\end{gather*}
Proposition \ref{i2j3} and our assumptions give $\lie d(V(\bomega''), V(\bomega_{i_1,j_1}))=0$. Since $\bos_{0,\varepsilon_s}(-1,\ell_s+1) = [i_1,j_1]\vee \bos_{0,\varepsilon_s}(0,\ell_s+1)$ the first equality follows from \eqref{omega''} and Proposition \ref{kkop}(iii).
Assuming  that the second equality holds we see that together with  Proposition \ref{kkop}(iii),(v) and Proposition \ref{i2j3} we get  \begin{gather*}\lie d(V(\bomega_{\bos_{0,\varepsilon_s}(-1, \ell_s+1)}), V(\bomega_{\bos_{0,\varepsilon_{s'}}(-1, \ell_{s'}+1)}))\le \lie d( V(\bomega_{\bos_{0,\varepsilon_{s}}(-1, \ell_s+1)}), V(\bomega_{i_1,j_1}))\\ \le 
\lie d(V(\bomega_{i_1,j_1}\bomega_{i_2,j_2}), V(\bomega_{i_1,j_1}))=0,
\end{gather*} and the third equality is proved.\\\\
The second equality is clearly a consequence of the following stronger proposition.
 \begin{prop} \label{irredcond} Assume that $r\ge 3$ and that        $2\le \ell_1,\, \ell_2\le r$ and  $\varepsilon_1,\varepsilon_2\in\{0,1\}$  are such that one of the following hold: \begin{gather*}
     2\ell_1+\varepsilon_1\ge  2\ell_2+\varepsilon_2\ \ {\rm and}\ \  \ell_2+\varepsilon_2\ \ {\rm is\ odd}\ \ {\rm or}\ \  2\ell_1+\varepsilon_1\le 2\ell_2+\varepsilon_2\ \ {\rm and}\ \ \ell_1+\varepsilon_1\ {\rm is\ even}.\end{gather*}
        Then $$
        \lie d(V(\bomega_{\bos_{0,\varepsilon_1}(-1, \ell_1+1)}), V(\bomega_{\bos_{0,\varepsilon_2}(0, \ell_2+1)}))=0.$$
        \end{prop}

        The proposition is proved in the rest of the section by induction on $\min\{\ell_1,\ell_2\}\ge 2 $.
        \subsubsection{} \label{first} The following result shows that the proposition holds when $\ell_1=2$ and $\varepsilon_1=0$. \begin{lem} Suppose that $r\ge 3$. For $\varepsilon\in\{0,1\}$ and $2\le \ell\le r$ with $\ell+\varepsilon\ge 3$ we have  $$\lie d(V(\bomega_{i_1,j_1}\bomega_{i_2,j_2}), V(\bomega_{\bos_{0,\varepsilon}(0,\ell+1)}))=0.$$\end{lem}
    \begin{pf}
    If $\ell+\varepsilon=3$ then $(\ell, \varepsilon)\in\{(3,0), (2,1)\}$ and we have to  prove that \begin{equation}\label{case1}\lie d(V(\bomega_{i_1,j_1}\bomega_{i_2,j_2}), V(\bomega_{i_2,j_2}\bomega_{i_3,j_3}))=0= \lie d(V(\bomega_{i_1,j_1}\bomega_{i_2,j_2}), V(\bomega_{i_2,j_2}\bomega_{i_4,j_3})).\end{equation} Here we understand that the second equality  exists only if  $r\ge 4$.
   Suppose first that  $\epsilon_1=0$ and $\ell=3$ then the  discussion in Section \ref{weyl} shows that the     modules $$V(\bomega_{i_p,j_p})\otimes V(\bomega_{i_s,j_s}\bomega_{i_2, j_2})\otimes V(\bomega_{i_2,j_2}),\ \ \ \{s,p\}=\{1,3\}$$ are $\ell$--highest weight.  Moreover, the modules $V(\bomega_{i_s,j_s}\bomega_{i_2, j_2})\otimes V(\bomega_{i_2,j_2})$ are irreducible by Proposition \ref{kkop}(v) for $s=1,3$. It follows that the modules $V(\bomega_{i_s, j_s}\bomega_{i_2,j_2})\otimes V(\bomega_{i_p,j_p}\bomega_{i_2,j_2})$ are $\ell$--highest weight for $\{s,p\}=\{1,3\}$ and hence irreducible. 
    If $\epsilon_1=1$ and $(\ell,\varepsilon)=(3,0)$
   then the module $$V(\bomega_{i_2,j_2})\otimes V(\bomega_{i_2,j_2}\bomega_{i_s,j_s})\otimes V(\bomega_{i_p,j_p}),\ \ \{s,p\}=\{1,3\}$$ is $\ell$--highest weight and the argument is completed as before. \\\\
   If $(\ell, \varepsilon)=(2,1)$  then 
        the arguments  are similar: if $\epsilon_1=0$ then  one works with   the $\ell$--highest weight module $V(\bomega_{i,j})\otimes V (\bomega_{i_2,j_2}\bomega_{i',j'})\otimes V(\bomega_{i_2,j_2})$ with $\{[i,j], [i',j']\}=\{[i_1,j_1], [i_4, j_3]\}$  and if $\epsilon_1=1$ with $V(\bomega_{i_2,j_2})\otimes V (\bomega_{i_2,j_2}\bomega_{i,j})\otimes V (\bomega\bomega_{i',j'})$.
        \\\\
         If $\ell+\varepsilon>3$ then either $\ell>3$ or $\ell=3$ and $\varepsilon=1$. In this case we note  that $[i_s,j_s]$ for $s=1,2$ is not connected to $[i_p,j_\ell]$ with $p,\ell \ge 4$ and so $$\lie d(V(\bomega_{i_1,j_1}\bomega_{i_2,j_2}), V(\bomega_{\bos_{0,\varepsilon}(2,\ell+1)}))=0.$$
        Since $\bos_{0,\varepsilon}(0,\ell+1)=\bos(1,3)\vee\bos_{0,\varepsilon}(2,\ell+1)$ we use  the $\ell+\varepsilon=3$ case and  apply  Proposition \ref{kkop}(iii) to deduce the lemma in this case.
    \end{pf}
    \subsubsection{} We consider the case when $\ell_1=2$ and $\varepsilon_1=1$. Since $\ell_1+\varepsilon_1$ is odd we must have $2\ell_2+\varepsilon_2\le 5$. Since $\ell_2\ge 2$ we get  $\ell_2=2$ and $\varepsilon_2=1$ and we have to prove that
$$\lie d(V(\bomega_{i_1,j_1}\bomega_{i_2,j_2}\bomega_{i_4,j_3}),  V(\bomega_{i_2,j_2}\bomega_{i_4,j_3}))=0.$$ But this follows from \eqref{case1} and an application of Proposition \ref{kkop}(v).
    \subsubsection{} We show that the proposition holds when $\ell_2=2$ and $\ell_1\ge 3$. Since $2\ell_1+\varepsilon_1\ge 6\ge 4+\varepsilon_2$  
    we must have $2+\varepsilon_2$ is odd forcing $\varepsilon_2=1$. Recalling that,
    \begin{gather*}\bos_{0,1}(0,3)=([i_2,j_2], [i_4,j_3]),\ \ \bos_{0,\varepsilon_1}(-1,\ell_1+1)=([i_1,j_1],[i_2,j_2])\vee \bos_{0,\varepsilon_1}(1,\ell_1+1),\end{gather*}
    the following equalities establish this case:
    \begin{gather*}
    \lie d(V(\bomega_{i_2,j_2}\bomega_{i_4,j_3}), V(\bomega_{i_1,j_1}\bomega_{i_2,j_2}))=0=\lie d(  V(\bomega_{i_2,j_2}\bomega_{i_4,j_3}), V(\bomega_{\bos_{0,\varepsilon_1}(1,\ell_1+1)})).
    \end{gather*} Here, the first equality immediate from Lemma \ref{first} while  the second follows by an application of Proposition \ref{i2j3}(ii) to $\bos(1,r)$.
 \subsubsection{}\label{equal}
We prove the inductive step; assume that $\min\{\ell_1,\ell_2\}\ge 3$ and set
\begin{gather*}\tilde\bos=\bos(1,r),\ \ \tilde{\ell_1}=\ell_2-1,\ \ \tilde{\ell_2}=\ell_1-1. \end{gather*} It is trivial to check that  $$\tilde{\bos}_{0,\varepsilon_2}(-1,\tilde\ell_1+1)= \bos_{0,\varepsilon_2}(0,\ell_2+1),\ \ \ \ \tilde{\bos}_{0,\varepsilon_1}(0,\tilde\ell_2+1)=\bos_{0,\varepsilon_1}(1,\ell_1+1). $$
Moreover,
\begin{gather*}
2\ell_1+\varepsilon_1\ge 2\ell_2+\varepsilon_2 \ \ {\rm and}\ \ \ell_2+\varepsilon_2\ \ {\rm odd}\iff 2\tilde \ell_1+\varepsilon_2\le 2\tilde\ell_2+\varepsilon_1 \ \ {\rm and}\ \ \tilde\ell_1+\varepsilon_2\ \ {\rm even},\\
2\ell_1+\varepsilon_1\le 2\ell_2+\varepsilon_2 \ \ {\rm and}\ \ \ell_1+\varepsilon_1\ \ {\rm even}\iff 2\tilde \ell_1+\varepsilon_2\ge 2\tilde\ell_2+\varepsilon_1 \ \ {\rm and}\ \ \tilde\ell_2+\varepsilon_1\ \ {\rm odd}.
    \end{gather*}
Since  $\min\{\tilde\ell_1,\tilde\ell_2\}=\min\{\ell_1,\ell_2\}-1$, the inductive hypothesis applies to $\tilde\bos$ and so
\begin{equation}\label{crux}\lie d(V(\bomega_{\bos_{0,\varepsilon_1}(1,\ell_1+1)}), V(\bomega_{\bos_{0,\varepsilon_2}(0,\ell_2+1)}))=0.\end{equation}
By Lemma \ref{first}$$\lie d(V(\bomega_{i_1,j_1}\bomega_{i_2,j_2}), V(\bomega_{\bos_{0,\varepsilon_2}(0,\ell_2+1)}))=0$$  and so an application of Proposition \ref{kkop}(iii) gives the inductive step.

    \section{An isomorphism of Grothendieck rings}\label{isoGr}
    Let  $\bos\in\uS^s\cap \mathbb I_n^r$ and $\bos'\in\uS^s\cap\mathbb I_N^r$ for some $n, N\in\mathbb N$. We give  a necessary condition for the existence of a ring isomorphism between $\cal K(\mathscr C_n(\bos))$ and $\cal K(\mathscr C_N(\bos'))$ which maps the class of an irreducible representation   to the class of an irreducible representation.
    \subsection{} We state the  main result of this section.
    \begin{thm}\label{main} Suppose that $\bos$ and $\bos'$ are prime elements of $\uS^s$ such that
    \begin{gather} 
 \label{isoc0} \epsilon_m=\epsilon_m',\ \ 1\le m\le r,\\
 \label{isoc2}
         i_{m-1}=i_{m+2}\iff i'_{m-1}=i'_{m+2},\ \ j_{m-1}=j_{m+2}\iff j'_{m-1}=j'_{m+2},\ \ 2\le m\le r-2,
         \\ \label{isoc1} [i_m,j_\ell]\in\mathbb I_n(\bos)\iff [i_m',j_\ell']\in\mathbb I_N(\bos').\end{gather}
         \begin{enumerit}
             \item[(i)] The assignment $[i_m, j_\ell]\to [i_m',j_\ell']$  defines a bijection $\mathbb I_n(\bos)\to \mathbb I_N(\bos')$ and induces a bijection $\mathbb I_n^r(\bos)\to\mathbb I_N^r(\bos')$ such that $$\bop\bor(\bos)\sqcup\bof\bor(\bos)\leftrightarrow  \bop\bor(\bos')\sqcup\bof\bor(\bos').$$
   \item[(ii)] The assignment $\bomega_{i_m,j_\ell}\to \bomega_{i_m',j_\ell'}$ induces          an isomorphism of monoids $\eta_{\bos,\bos'}:\cal I_n^+(\bos)\to\cal I_N^+(\bos')$ and   an isomorphism of rings $\tilde{\eta}_{\bos,\bos'}: \cal K(\mathscr C_n(\bos))\to\cal K(\mathscr C_N(\bos'))$ such that
     \begin{gather} \label{irredimage}\tilde\eta_{\bos,\bos'}([V(\bomega)])= [V(\eta_{\bos, \bos'}(\bomega))],\ \ \bomega\in\cal I_n^+(\bos).\end{gather}
     \end{enumerit}
\end{thm}
\subsection{} We prove that the map $[i_m, j_\ell]\to[i_m', j_\ell']$ is well-defined. Suppose that $[i_m,j_\ell]\in\mathbb I_n(\bos)$ and  $i_m=i_p$ and $j_\ell=j_s$ for some $1\le p,s\le r$. By Lemma \ref{enumerate} it follows that either $\epsilon_m=0$ and $p=m+3$ or $\epsilon_p=0$ and $m=p+3$ and either  $\epsilon_\ell=1$ and $s=\ell+3$ or $\epsilon_s=1$ and $\ell=s+3$. It follows from \eqref{isoc2} that $i_m'=i_p'$ and $j_\ell'=j_s'$ as needed. The proof that the map is one--one is identical and \eqref{isoc1} shows that the map is onto.\\\\
Using \eqref{isoc0}-\eqref{isoc1} we see that the induced bijection $\mathbb I_n(\bos)\to\mathbb I_N(\bos')$ maps $\bop\bor(\bos)\sqcup\bof\bor(\bos)\mapsto \bop\bor(\bos')\sqcup\bof\bor(\bos')$ and the proof of part (i) is complete.\\\\
For part (ii), we note that the existence of the isomorphism $\eta_{\bos,\bos'}:\cal I_n^+(\bos)\to\cal I_N^+(\bos')$ is now immediate since the monoids are generated by $\bomega_{i_m,j_\ell}$ and $\bomega_{i_m',j_\ell'}$ with $[i_m, j_\ell]\in\mathbb I_n(\bos)$ and $[i_m, j_\ell]\in\mathbb I_N(\bos')$ respectively.\\\\
By \cite{FR99}, we know that  $\cal K(\mathscr C_n(\bos))$ is a polynomial ring in the variables $[V(\bomega_{i_m, j_\ell})]$, $[i_m,j_\ell]\in\mathbb I_n(\bos)$ and a similar statement for $\cal K(\mathscr C_N(\bos'))$. It  follows form part (i) of theorem that the assignment $[V(\bomega_{i_m, j_\ell})]\to [V(\bomega_{i_m', j_\ell'})]$  defines   an isomorphism $\tilde{\eta}_{\bos,\bos'}: \cal K(\mathscr C_n(\bos))\to \cal K(\mathscr C_N(\bos'))$ of rings.\\\\
The rest of the section is devoted to proving that $\tilde{\eta}_{\bos,\bos'}$ satisfies \eqref{irredimage}; i.e., that it takes the class of an irreducible module to the class of an irreducible module.
\subsection{} \label{etarest} We make some preliminary observations. Note that the elements $\bos_{\varepsilon,\varepsilon'}(p,\ell+1)$ and $\bos'_{\varepsilon,\varepsilon'}(p,\ell+1)$ also satisfy \eqref{isoc0}-\eqref{isoc1}. Hence we have the corresponding isomorphisms
\begin{gather*}\eta_{\bos_{\varepsilon,\varepsilon'}(p,\ell+1),  \bos'_{\varepsilon,\varepsilon'}(p,\ell+1)}:\cal I_n^+(\bos_{\varepsilon,\varepsilon'}(p,\ell+1))\to \cal I_N^+(\bos'_{\varepsilon,\varepsilon'}(p,\ell+1)),\\ \tilde{\eta}_{\bos_{\varepsilon,\varepsilon'}(p,\ell+1), \bos'_{\varepsilon,\varepsilon'}(p,\ell+1)}:\cal K(\mathscr C_n(\bos_{\varepsilon,\varepsilon'}(p,\ell+1)))\to\cal K(\mathscr C_N(\bos'_{\varepsilon,\varepsilon'}(p,\ell+1)).\end{gather*}
Since 
$$\mathbb I_n(\bos_{\varepsilon,\varepsilon'}(p,\ell+1))\subset \mathbb I_n(\bos),\ \ \mathbb I_N(\bos'_{\varepsilon,\varepsilon'}(p,\ell+1))\subset \mathbb I_N(\bos'),$$ it is now clear that the restriction of $\tilde\eta_{\bos,\bos'}$ to $\cal K(\mathscr C_n(\bos_{\varepsilon,\varepsilon'}(p,\ell+1)))$ is $\tilde{\eta}_{\bos_{\varepsilon,\varepsilon'}(p,\ell+1), \bos'_{\varepsilon,\varepsilon'}(p,\ell+1)}$. A similar discussion applies if we work with elements of $\bof\bor(\bos)$.
\subsection{}  The first and second identities in the next proposition follow from  \cite[Proposition 9]{BC25a} by taking  $r_1=2$, and the third one follows from \cite[Lemma 8.3]{BC25a}. Recall that $\epsilon_1\in\{0,1\}$ is chosen so that $i_{1+\epsilon_1}>i_{2-\epsilon_1}$. In the following we shall assume that $\bos(2,r)$ appears only if $r\ge 3$.
\begin{prop}\label{bcss}
Suppose that $\bos\in\uS^s$ is prime. Then 
\begin{gather}\label{clust1}[V(\bomega_{i_1,j_1})][V(\bomega_{\bos(1,r)})]= [V(\bomega_{\bos})]+ [V(\bomega_{i_2,j_1}\bomega_{i_1,j_2}\bomega_{\bos(2,r)})],\\
 \label{clust2}[V(\bomega_{i_2,j_1}\bomega_{i_1,j_2}\bomega_{\bos(2,r)})]= [V(\bomega_{i_2,j_1}^{1-\epsilon_1}\bomega_{i_1,j_2}^{\epsilon_1})] [V(\bomega_{i_2,j_1}^{\epsilon_1}\bomega_{i_1,j_2}^{1-\epsilon_1} \bomega_{\bos(2,r)})].\end{gather} If in addition $i_1=i_4$ or $j_1=j_4$ then 
 \begin{gather}\label{clust3} [V(\bomega_{i_2,j_1}^{\epsilon_1}\bomega_{i_1,j_2}^{1-\epsilon_1} \bomega_{\bos(2,r)})]  = [V(\bomega_{i_2,j_1}^{\epsilon_1}\bomega_{i_1,j_2}^{1-\epsilon_1}\bomega_{i_3,j_3})][V(\bomega_{\bos(3,r)})].\end{gather}
\hfill\qedsymbol
\end{prop}
\subsection{} We prove that \eqref{irredimage} holds by induction on $r$ with induction obviously beginning when $r=1$. If $r=2$ then 
using
\eqref{clust1}-\eqref{clust2} we have 
\begin{gather*} [V(\bomega_{i_1,j_1})][V(\bomega_{i_2,j_2})]=[V(\bomega_{i_1,j_1}\bomega_{i_2,j_2})]+[V(\bomega_{i_1,j_2})][V(\bomega_{i_2,j_1})],\\
[V(\bomega_{i_1',j_1'})][V(\bomega_{i_2',j_2'})]=[V(\bomega_{i_1',j_1'}\bomega_{i_2',j_2'})]+[V(\bomega_{i_1',j_2'})][V(\bomega_{i_2',j_1'})].
\end{gather*}
Applying $\tilde\eta$ to both sides of the first equation and using the result for $r=1$ proves that $\tilde\eta([V(\bomega_{i_1,j_1}\bomega_{i_2,j_2})])=[V(\bomega_{i_1',j_1'}\bomega_{i_2',j_2'})]$. 
\\\\
If $\bomega\in\cal I_n^+(\bos)$ then writing $$\bomega= \bomega_{i_1,j_1}^{a_1}(\bomega_{i_1,j_1}\bomega_{i_2,j_2})^a\bomega_{i_2,j_2}^{a_2}\bomega_{i_1,j_2}^b\bomega_{i_2, j_1}^c,\ \ {\rm where}\ \  a_1=0\ \ {\rm or} \ \ a_2=0,$$
it follows from  Proposition \ref{kkop}(iv) and the fact that $i_1<i_2<j_1<j_2$ or $i_2<i_1\le j_2<j_1$ that $$[V(\bomega)]=[V(\bomega_{i_1,j_1})]^{a_1}[V(\bomega_{i_1,j_1}\bomega_{i_2,j_2})]^a[V(\bomega_{i_2,j_2})]^{a_2}[V(\bomega_{i_1,j_2})]^b[V(\bomega_{i_2, j_1})]^c.$$ Since  an identical factorization holds for $V(\eta(\bomega))$ the  proof in the  case $r=2$ follows. In particular, using the discussion in Section \ref{etarest}, equation  \eqref{irredimage} holds if $\bomega=\bomega_\bos$ with $\bos\in\bof\bor(\bos)$.

\subsection{}  
Assume that the  result holds for all $\bomega'\in\cal I_n^+(\bos_1)$ with $\bos_1\in \uS^s\cap\mathbb I_n^{r-1}$ prime. \\\\
Suppose first that $\bomega=\bomega_{\bos_{\varepsilon,\varepsilon'}(p,\ell+1)}$ where either $p>-1$ or $\ell <r$.
The inductive hypothesis applies to $\bomega$ and the   discussion in Section \ref{etarest} shows that $$\tilde{\eta}_{\bos,\bos'}([V(\bomega)])=\tilde\eta_{\bos_{\varepsilon, \varepsilon'}(p,\ell+1)}
([V(\bomega)])=[V(\eta_{\bos'_{\varepsilon, \varepsilon'}(p,\ell+1)}(\bomega))]=[V(\eta_{\bos,\bos'}(\bomega))],$$ as needed.
Hence it suffices to prove that $\tilde\eta_{\bos,\bos'}([V(\bomega_\bos)])=[V(\bomega_{\bos'})]$.
If  $i_1=i_4$ or $j_1=j_4$ then  the inductive hypothesis applies to \eqref{clust3} and gives
$$\tilde\eta_{\bos,\bos'}([V(\bomega_{i_2,j_1}^{\epsilon_1}\bomega_{i_1,j_2}^{1-\epsilon_1} \bomega_{\bos(2,r)})])  = [V(\bomega_{i_2',j_1'}^{\epsilon_1'}\bomega_{i_1',j_2'}^{1-\epsilon_1'}\bomega_{i_3',j_3'})][V(\bomega_{\bos'(3,r)})].$$ Otherwise, \eqref{clust2} and gives,
$$\tilde\eta_{\bos,\bos'}([V(\bomega_{i_2,j_1}\bomega_{i_1,j_2} \bomega_{\bos(2,r)})])=[V(\bomega_{i_2',j_1'}^{1-\epsilon_1'}\bomega_{i_1',j_2'}^{\epsilon_1'})] [V(\bomega_{i_2',j_1'}^{\epsilon_1'}\bomega_{i_1',j_2'}^{1-\epsilon_1'} \bomega_{\bos'(2,r)})].$$
Applying $\tilde\eta_{\bos,\bos'}$ to both sides of \eqref{clust1} and using  the preceding observations gives \eqref{irredimage} for $\bomega_\bos$.\\\\
Now assume that $\bomega$ is not necessarily prime and also $\bomega\notin \cal I_n^+(\bos(1,r))$. We proceed by induction on $\Ht\bomega$.
The results of Section \ref{allprime} and the definition of $\eta_{\bos,\bos'}$  show that  if $$(*)\ \  \ \ \bomega\bomega_1^{-1}\in\cal I_n^+,\ \ \bomega_1\in\{\bomega_{i_2,j_1},\ \ \bomega_{i_1,j_3}\bomega_{i_2,j_2}, \ \ \bomega_{i_1,j_1}\bomega_{i_2,j_3}\}$$ we have $$V(\bomega)\cong V(\bomega_1)\otimes V(\bomega\bomega_1^{-1}),\ \ V({\eta}_{\bos,\bos'}(\bomega))\cong V({\eta}_{\bos,\bos'}(\bomega_1)))\otimes V({\eta}_{\bos,\bos'}(\bomega\bomega_1^{-1})).$$ Since we have proved the result for elements in $\bop\bor(\bos)$ and $\Ht\bomega\bomega_1^{-1}<\Ht\bomega$ we have
$$\tilde {\eta}_{\bos,\bos'}([V(\bomega)])= \tilde {\eta}_{\bos,\bos'}([V(\bomega_1)])\tilde {\eta}_{\bos,\bos'}([V(\bomega\bomega_1^{-1})])=[V({\eta}_{\bos,\bos'}(\bomega_1))][V({\eta}_{\bos,\bos'}(\bomega\bomega_1^{-1}))] = [V({\eta}_{\bos,\bos'}(\bomega))].$$
If $\bomega$ does not satisfy $(*)$ but does satisfy $$(**)\ \ \bomega\bomega_1^{-1}\in\cal I_n^+,\ \ \bomega_1\in\{\bomega_{i_1,j_3}\ \bomega_{i_2,j_3}\}$$
then an identical argument establishes \eqref{irredimage}.
\\\\
Otherwise, we are in situation of Section \ref{i1j2} and have  \begin{equation*}
\bomega=\bomega_{i_1,j_1}^{a_1}\bomega_{i_1,j_2}^{a_2}\bomega_{i_2, j_2}^c\bomega_2,\ \ \bomega_2\in \cal I_n^+(\bos_{1,0}(1,r+1)),\ \ c>0.
\end{equation*} and a similar statement for $\eta_{\bos,\bos'}(\bomega)$. \\

Assuming that $a_2>0$ and letting $0\leq d,d'\leq a_2$, $s_k\in\mathbb N$, $\varepsilon_p$, $1\leq k\leq a_2-d-d'$ and $\bomega_3\bomega_{i_1,j_2}^{-1}\notin\cal I_n^+$ as in \eqref{a2>0case} we have
\begin{equation}\label{a2>01}
    [V(\bomega)] =[V(\bomega_{i_1,j_2})]^{d} [V(\bomega_{i_1,j_2}\bomega_{i_4, j_3})]^{d'}\left(\prod_{k=1}^{a_2-d-d'} [V(\bomega_{\bos_{1,\varepsilon_k}(1,s_k+1)})]\right)[V(\bomega_{i_1,j_1}^{a_1}\bomega_{i_2,j_2}^c\bomega_3)],
\end{equation}
Similarly, 
\begin{equation}\label{a2>02}[V(\eta_{\bos,\bos'}(\bomega))] =[V(\bomega_{i_1',j_2'})]^{d} [V(\bomega_{i_1',j_2'}\bomega_{i_4', j_3'})]^{d'}\left(\prod_{k=1}^{a_2-d-d'} [V(\bomega_{\bos'_{1,\varepsilon_k}(1,s_k+1)})]\right)[V(\bomega_{i_1',j_1'}^{a_1}\bomega_{i_2',j_2'}^c\eta_{\bos,\bos'}(\bomega_3))].\end{equation}
Applying $\tilde\eta_{\bos,\bos'}$ to both sides of the \eqref{a2>01} and using the inductive hypothesis gives that its right hand side is equal to the right hand side of \eqref{a2>02} and hence we get $\tilde\eta_{\bos,\bos'}([V(\bomega)])=[V(\eta_{\bos,\bos'}(\bomega))]$, which establishes the inductive step in this case. 
\\\\
If $a_2=0$ we are finally we are in the situation of Section \ref{finalcase}; the argument here is identical to the one given above and we omit the details.

\section{Monoidal Categorification}\label{moncat}
In this section we  prove that if $\bos\in\uS^s$ is prime then $\cal K(\mathscr C_n(\bos))$ is the monoidal categorification of a cluster algebra. Throughout this section for $m\ge 1$ we set $I_m=\{1,2,\cdots, m\}$.
\subsection{A height function} Given $\bos=([i_1,j_1],\cdots, [i_r, j_r])\in\uS^s$ prime and $r\ge 3$, set
$$N(\bos)=r+\#\{2\le s\le r-2: i_{s-1}\ne i_{s+2}\ \ {\rm and}\ \ j_{s-1}\ne j_{s+2}\}.$$ For ease of notation and when there is no ambiguity we set $N=N(\bos)$.
Define a subset $\{p_1,\cdots, p_r\}$ of $\mathbb Z$ setting
\begin{gather*}
    p_1=1,\ \  p_{2}=2,\\ 
p_{m+1}=p_{m}+\begin{cases} 2,\ \ i_{r-m+2}\ne i_{r-m-1}\ \ {\rm and}\ \ j_{r-m+2}\ne j_{r-m-1},\\
    1,\ \ {\rm otherwise},
        \end{cases}\ \ 2\le m\le r-2,\\
        \ p_r=p_{r-1}+1.
\end{gather*}
\begin{lem}\label{pi}
    We have $p_r=N$ or equivalently  $p_{r-1}=N-1$. In particular, $\{p_1,\cdots, p_r\}\subset I_N$.
\end{lem}
\begin{pf} The lemma is proved by induction on $r$ with the case $r=3$ being clear. Assuming the lemma  for $\bos'=\bos(0,r-1)$ we note that $$N(\bos)= \begin{cases} N(\bos')+2,\ \ i_{r-3}\ne i_r \ \ {\rm and }\ \ j_{r-3}\ne j_r,\\
N(\bos')+1,\ \ {\rm otherwise}.
\end{cases}$$
Moreover if we let $\{p_1',\cdots, p_{r-1}'\}$ be the subset corresponding to $\bos'$ we have by the inductive hypothesis that $p_s'=p_s$ for $1\le s\le r-2$ and $p_{r-2}'=N(\bos')-1$.     It follows now from the  definition $p_{r-1}=N(\bos)-1$ as needed and the proof of the lemma is complete.
\end{pf}

\subsection{The category $\mathscr C(\xi_\bos)$} Recall that a function $\xi: I_M\to\mathbb Z $ is called a height function if $|\xi(s)-\xi(s+1)|\le 1$ for all $1\le s\le M-1$. Recall also from Definition \ref{altsnakedef} the elements $\{\epsilon_s : 1\leq s\leq r\}$ associated to an element $\bos\in \uS^s$. For $\bos\in\uS^s$ prime, set $N=N(\bos)$; it is easy to check that the following formulae define a height function $\xi_\bos: I_{N}\to\mathbb Z$:
\begin{gather}\label{defxis}
    \xi_\bos(p_r)=p_r,\ \ \xi_\bos(p_m-1)=\xi_\bos(p_m+1),\ \ \xi_\bos(p_{m+1-\epsilon_{r-m}})-\xi_\bos(p_{m+\epsilon_{r-m}})=p_{m+1}-p_m.
\end{gather} 

\iffalse\begin{ex}
    If $r=4$ and $a_1\neq a_4$ for $a\in\{i,j\}$ then 
    $$N(\bos)=5, \ \ p_1=1, \ \ p_2=2, \ \ p_3=4, \ \ p_4=5.$$
    Moreover, if $\epsilon_1=0$ then we have $\epsilon_1=\epsilon_3=0$ and $\epsilon_2=\epsilon_4=1$ and
    $$\xi_\bos(5)=5, \ \ \xi_\bos(4)= 4, \ \ \xi_\bos(3)=5, \ \ \xi_\bos(2)=6, \ \ \xi_\bos(1)=5.$$
\end{ex}\fi 

For $1\le s\le N$  notice that $\xi_\bos(s)-s\in 2\mathbb Z$ and set
  \begin{gather*}
      2i^{\xi_\bos}_s=\xi_\bos(s)-s, \ \ 2j^{\xi_\bos}_s=\xi_\bos(s)+s.
      \end{gather*}
 Set 
  $$\mathbb I_{N}(\xi_\bos)=\{[i_s^{\xi_\bos},\   j_s^{\xi_\bos}], [i_s^{\xi_\bos}-1,\   j_s^{\xi_\bos}-1]: 1\le s\le N\},$$
 and let  $\cal I_{N}^+(\xi_\bos)$ be the submonoid of $\cal I_{N}^+$ generated by elements  $\bomega_{i,j}$, $[i,j]\in\mathbb I_{N}(\xi_\bos)$. \\

\begin{ex}
    Suppose that $r=5$, and $\bos\in\uS^s\cap \mathbb I_n^5$ is prime such that $\epsilon_1=0$, $a_1\neq a_4$, $a\in \{i,j\}$ and $j_2=j_5$. Then 
    \begin{gather*}
        N(\bos)=6, \ \ p_1=1,\ \ p_2=2, \ \ p_3= 3, \ \ p_4= 5,\ \ p_5= 6,\\
        \xi_\bos(6)= 6 ,\ \ \xi_\bos(5)= 5 ,\ \ \xi_\bos(4)= 6,\ \ \xi_\bos(3)= 7,\ \ \xi_\bos(2)= 6 ,\ \ \xi_\bos(1)= 7,\\
        \mathbb I_{N(\bos)}(\xi_\bos) = \{[3,4], [2,3], [2,4], [1,3], [2,5],[1,4],[1,5], [0,4], [0,5], [-1,4], [0,6],[-1,5]\}.
    \end{gather*}
\end{ex}

 In \cite{HL10} the authors defined (more generally for an arbitrary height function) the subcategory $\mathscr C_{\xi_\bos}$ of $\mathscr F_{N}$ to be  the full subcategory of modules whose Jordan--H\"older constituents are of the form $V(\bomega)$ with $\bomega\in\cal I_{N}^+(\xi_\bos)$. The following was proved in the language of Drinfeld polynomials in \cite{HL10} for the bipartite height function and in \cite{BC19a} in general. 
\begin{thm}\label{moncatht}
  \begin{enumerit}
      \item [(i)] The category $\mathscr C_{\xi_\bos}$ is a monoidal tensor category and the ring $\cal K(\mathscr C_{\xi_\bos})$ is a monoidal categorification of a cluster algebra of type $A_{N(\bos)}$ with $N(\bos)$ frozen variables.
      \item[(ii)] The cluster variables (resp. cluster monomials) are in bijective correspondence with the prime objects of $\mathscr C_{\xi_\bos}$  (resp. irreducible objects of $\mathscr C_{\xi_\bos}$).
      \end{enumerit}
      \hfill\qedsymbol
      \end{thm}
      The index set for prime and  modules in $\mathscr C_{\xi_\bos}$ were determined in \cite{BC19a} in the language of Drinfeld polynomials. We now give the index set in the language of this paper. \\\\
     For $1\le m\le \ell\le r$ set 
     $$\bomega(p_m,p_\ell) =\prod_{k=m}^\ell \bomega_{i_{p_k}^\xi -\epsilon_{r-k+1},\, j_{p_k}^\xi -\epsilon_{r-k+1}}$$
     Given $1\le s<s'\le N(\bos)$ let $1\leq m,\ell\leq r$ be the unique elements satisfying
     \begin{equation}\label{ss'bounds}1\leq m\leq \ell\leq r \ \ {\rm and} \ \ p_{m-1}<s\leq p_m\leq p_\ell\leq s'<p_{\ell+1},
     \end{equation}
     and set
\begin{gather*}\bomega(s,s')=\bomega_{i^\xi_{p_{m-1-\epsilon_{r-m+1}}},\, j^\xi_{p_{m-2+\epsilon_{r-m+1}}}}^{1-\delta_{s,p_m}} \bomega(p_m, p_\ell)\ \ \bomega_{i^\xi_{p_{\ell +2-\epsilon_{r-\ell+1}}},\, j^\xi_{p_{\ell+1+\epsilon_{r-\ell+1}}}}^{1-\delta_{s',p_\ell}}.
     \end{gather*}
    Then, the index set $\bp\bor(\xi_\bos)$ is the union of the following two sets:
    \begin{gather*}
      \{\bomega_{i^\xi_s- 1, j_s^\xi- 1},  \ \bomega_{i^\xi_s, j_s^\xi}: \ 1\le s\le N(\bos)\}\cup 
      \{\bomega(s,s'): 1\le s<s'\le N(\bos)\}, 
 \end{gather*} 
 and the frozen elements $$\bold F\bor(\xi_\bos)=\{\bomega_{i_s^\xi, j_s^\xi} \bomega_{i_s^\xi -1, j_s^\xi-1}: 1\le s\le N(\bos)\}.$$
 
\subsection{} Part (iv) of Theorem \ref{mainthm} is immediate from the following proposition.
      \begin{prop}\label{red}
           Suppose that $\bos=([i_1,j_1],\cdots, [i_r,j_r])\in\uS^s\cap \mathbb I_n^r$ is prime and stable with $j_{\min}=i_{\max}$ and $n+1=j_{\max}-i_{\min}$. 
           There exists an isomorphism $\eta_{\xi, \bos}:\cal I_{N(\bos)}^+(\xi_\bos)\to\cal I_n^+(\bos)$ of monoids and a corresponding isomorphism of rings $\tilde\eta_{\xi, \bos}: \cal K(\mathscr  C_{\xi_\bos})\to\cal K(\mathscr C_n(\bos))$ which maps irreducible (prime) objects onto irreducible (prime) objects. 
       \end{prop} 
       The proof of the proposition is given in the rest of the section.
\begin{rem} Note that for $1\leq s<s'\leq N(\bos)$ and $1\leq m\leq \ell\leq r$ as in \eqref{ss'bounds} the assignment $$\bomega(s,s')\mapsto \bomega_{\tilde\bos}, \ \ \tilde\bos = \bos_{1-\delta_{s',p_\ell},1-\delta_{s,p_m}}(p_{r-\ell+1}-2, p_{r-m+1})$$
defines a bijection between $\bp\bor(\xi_\bos)\to \bp\bor(\bos)$\end{rem}

\subsubsection{}   We deduce the following consequence of Theorem \ref{moncatht}.
\begin{prop} Suppose that  $([i_m,j_\ell], [i_p,j_s])\in\mathbb I_{N(\bos)}(\xi_\bos)^2$ is connected with $i_p<i_m$. If $j_\ell-i_p\le N(\bos)$ (resp. $i_m\ne j_s$) then the interval $[i_p, j_\ell]$ (resp. $[i_m, j_s]$) is  in $\mathbb I_{N(\bos)}(\xi_\bos)$.
 \end{prop}
    \begin{pf}
 Working with  $([i_m, j_\ell], [i_p,j_s])\in\uS^s$, we have from Proposition \ref{bcss} that \begin{gather*}[V(\bomega_{i_m,j_\ell})][V(\bomega_{i_p,j_s})]= [V(\bomega_{i_m,j_\ell}\bomega_{i_p,j_s})]+ [V(\bomega_{i_m,j_s}\bomega_{i_p,j_\ell})],\\ [V(\bomega_{i_m,j_s}\bomega_{i_p,j_\ell})]=[V(\bomega_{i_m,j_s})][V(\bomega_{i_p,j_\ell})].\end{gather*} Since $\mathscr C_{\xi_\bos}$ is a monoidal tensor category the first equation shows that  the product $\bomega_{i_m,j_s} \bomega_{i_p,j_\ell}$ is an element of $\mathcal I_{N(\bos)}^+(\xi_\bos)$. On the other hand Theorem \ref{moncatht} implies that every irreducible representation in $\mathscr C_{\xi_\bos}$ is isomorphic to a tensor product of prime irreducible representations  in $\mathscr C_{\xi_\bos}$. The second equation shows that 
 the module $V(\bomega_{i_m, j_s}\bomega_{i_p,j_\ell})$ is not prime in $\mathscr F_{N(\bos)}$ and hence not prime in $\mathscr C_{\xi_\bos}$. It is now clear that the modules $V(\bomega_{i_m,j_s})$ and $V(\bomega_{i_p, j_\ell})$ must be objects of $\mathscr C_{\xi_\bos}$. In particular, 
we have 
       $$\bomega_{i_m,j_s}\neq \bold 1\implies i_m\ne j_s,\ \ \bomega_{i_p, j_\ell}\ne \bold 1\implies j_\ell-i_p\le n,
       $$
       and the proof of the proposition is complete. \end{pf}
\subsubsection{} Set
\begin{gather}\label{y'def}
s_m=p_{r-m+1},\ \ 
a_{s_m}'=a_{s_m}^{\xi_\bos}-\epsilon_{m},\ \ a\in\{i,j\},\ \ 1\le m\le r.\end{gather}
Using \eqref{defxis} we have 
\begin{gather}\label{distance}
    0<s_m-s_{m+1} = \xi_\bos(s_{m+\epsilon_m})-\xi_\bos(s_{m+1-\epsilon_m}), \ \  \ 1\leq m<r\\ \label{a=a}
     s_m-s_{m+1}=1\iff a_{m-1}= a_{m+2},\ \ {\rm for \ some}  \ a\in \{i,j\}, \ \   2\leq m\leq r-2.
\end{gather}
\begin{lem}\label{sxi}
    Retain the notation established so far. Then $$\bos_{\xi_\bos}:=([i_{s_1}', j_{s_1}'],\cdots,[i_{s_r}',j_{s_r}'])$$ is a prime element of $\uS^s$ and the pair $(\bos,\bos_{\xi_\bos})$ satisfy \eqref{isoc0}-\eqref{isoc1}.
  
\end{lem}
\begin{pf} Using \eqref{y'def} and \eqref{distance} the following inequalities are easily checked and prove that   $\bos_{\xi_\bos}\in\uS^s$ and is prime:
\begin{gather*} a'_{s_{m+\epsilon_m}}>a'_{s_{m+1-\epsilon_m}},\ \ \ a\in\{i,j\},\ \ 1\le p\le r-1,\\
     j_{s_m}'-i_{s_m}
     '\ne  j_{s_p}'-i_{s_p}',\ \ \ m\ne p,\\
    i'_{s_{m+1-\epsilon_{m}}}<i'_{s_{m+\epsilon_{m}}}\le j'_{s_{m+1-\epsilon_m}}<j'_{s_{m+\epsilon_{m}}}, \ \ 1\leq m\leq r-1\\
    i'_{s_m}<i'_{s_{m+2}}<j'_{s_{m+2}}<j_{s_m}',\ \ \ i'_{s_m}\le i_{s_{m+3}}'<j_{s_{m+3}}'\le j_{s_{m}}'.
\end{gather*}
Lemma \ref{enumerate}  gives
    $$j_{\max}'= j_{s_{1+\epsilon_1}}', \ \ j_{\min}'= j_{s_{r-1+\epsilon_r}}', \ \ i_{\min}'= i_{s_{2-\epsilon_1}}', \ \ i_{\max}'= i_{s_{r-\epsilon_r}}'.$$
Moreover, since $(s_1,s_2)=(N(\bos),N(\bos)-1)$ and $(s_{r-1}, s_r)=(2,1)$ one easily checks that
\begin{equation}\label{ijN}j'_{\min}=i'_{\max} \ \  {\rm and} \ \ j'_{\max}-i'_{\min}= N(\bos)+1.\end{equation} \\
Setting $\bos'=\bos_{\xi_\bos}$, it is clear by definition that $\epsilon_1=\epsilon_1'$ and hence \eqref{isoc0} holds.
Notice that \eqref{in} shows that the pair satisfies \eqref{isoc1}. 
Next we prove that  \eqref{isoc2} holds; using \eqref{distance} notice that
 \begin{gather}\label{distanceA}
   \xi(s_p)+s_p = \xi(s_{p-1})+s_{p-1} \ \ {\rm if} \ \ \epsilon_p=0, \ \ {\rm and} \ \ \xi(s_p)-s_p =\xi(s_{p-1})-s_{p-1} \ \ {\rm if} \ \ \epsilon_p=1.
 \end{gather}
Then, if $\epsilon_p=0$ we have 
\begin{gather*}
    j_{p-1}=j_{p+2} \overset{\eqref{a=a}}{\iff} s_p-s_{p+1}=1 \overset{\eqref{distance}}{\iff} \xi(s_p)-\xi(s_{p+1})=1\\ 
    \overset{\eqref{distanceA}}{\iff} \xi(s_{p+2})+ s_{p+2}+1 = \xi(s_{p-1})+s_{p-1}-1,\ \ {\rm i.e.}\ \ j'_{s_{p-1}}=j'_{s_{p+2}},\ \ \ 
\end{gather*}
as desired. If $\epsilon_p=1$ then
\begin{gather*}
    i_{p-1}=i_{p+1} \overset{\eqref{a=a}}{\iff} s_p-s_{p+1}=1 \overset{\eqref{distance}}{\iff} \xi(s_p)-\xi(s_{p+1})=-1\\
    \overset{\eqref{distanceA}}{\iff} \xi(s_{p+2})-s_{p+2}-1 = \xi(s_{p-1})-s_{p-1}+1\ \ {\rm i.e.}\ \ i'_{s_{p-1}}=i'_{s_{p+2}}
\end{gather*}
which completes the proof. \end{pf}

\begin{ex}
    If $\bos$ is as in Example \ref{defxis} then $\bos_{\xi_\bos} = ([0,6],[-1,4],[2,5],[1,3],[3,4])\in \uS^s\cap \mathbb I_6^5$ is prime.
\end{ex}

\subsubsection{} As a consequence of Lemma \ref{sxi} and Theorem \ref{main} we now have that $\cal K(\mathscr C_n(\bos))$ and $\cal K(\mathscr C_{N(\bos)}(\bos_{\xi_\bos}))$  are isomorphic via an isomorphism that maps prime objects to prime objects and irreducible modules to irreducible modules. Proposition \ref{red} obviously follows if we establish that the category $\mathscr C_{N(\bos)}(\bos_{\xi_\bos})$ and $\mathscr C_{\xi_\bos}$ are the same. This is true once we prove that 
\begin{equation}\label{sxis}\mathbb I_{N(\bos)}(\xi_\bos)=\mathbb I_{N(\bos)}(\bos_{\xi_\bos}).\end{equation}
We check first that $\mathbb I_N(\xi_\bos)\subset \mathbb I_N(\bos_{\xi_\bos})$. By definition we have
$[i_{s_p}^\xi-\epsilon_p ,j_{s_p}^\xi -\epsilon_p]=[i_{s_p}',\ j_{s_p}']$ and hence we have to prove the following two assertions:
\begin{gather}
    \label{turningpoints}[i_{s_p}^\xi+\epsilon_p-1 ,j_{s_p}^\xi +\epsilon_p-1]\in\mathbb I_N(\bos_{\xi_\bos}),\\
    \label{notturningpoints}s_p+1\ne s_{p-1}\implies   [i^\xi_{s_p+1}-\epsilon, j^\xi_{s_p+1}-\epsilon]\in\mathbb I_N(\bos_{\xi_\bos}),\ \ \epsilon\in\{0,1\}.
\end{gather}
A calculation using \eqref{distance}  shows that
\begin{gather*}
     [i_{s_1}^\xi+\epsilon_1-1,j_{s_1}^\xi+\epsilon_1-1] = [i_{s_{2+\epsilon_1}}',j_{s_{3-\epsilon_1}}'], \ \ \ [i_{s_r}^\xi+\epsilon_r-1,j_{s_r}^\xi+\epsilon_r-1] =[i_{s_{r-2+\epsilon_r}}', j_{s_{r-1-\epsilon_r}}']
\end{gather*}
and for $1<p<r$, 
\begin{gather*}
  [i_{s_p}^\xi -1, j_{s_p}^\xi-1] = [i_{s_{p+1}}', j_{s_{p-1}}']  \ \ {\rm if}\  \  \epsilon_p=0, \ \ {\rm and} \ \ 
[i_{s_p}^\xi, j_{s_p}^\xi] = [i_{s_{p-1}}', j_{s_{p+1}}'] \ \ {\rm if} \ \ \epsilon_p=1.
\end{gather*}
It is now easy to check that \eqref{turningpoints} holds. We turn to the proof of \eqref{notturningpoints}. Suppose that $s_p-1\neq s_p$, or equivalently that $s_{p-1}-s_p=2$; if $\epsilon_p=0$ then by \eqref{defxis} we have
$\xi(s_p)-\xi(s_{p}+1)=1$ and we get:
\begin{gather*}
    2i_{s_p+1}^\xi = \xi(s_p+1)-(s_p+1)=\xi(s_{p})-s_{p}-2\overset{\eqref{distanceA}}{=}\xi(s_{p+1})-s_{p+1}-2 = 2i_{s_{p+1}}',\\
    2(i_{s_p+1}^\xi-1)= \xi(s_p)-s_p-4 = \xi(s_{p-1})-s_{p-1} \overset{\eqref{distanceA}}{=} \xi(s_{p-2})-s_{p-2}=2i_{s_{p-2}}',\\
    2j_{s_p+1}^\xi = \xi(s_p+1)+(s_p+1)=\xi(s_p)+s_p =  2j_{s_{p}}',\\
2(j_{s_p+1}^\xi-1) = \xi(s_p+1)+(s_p+1)-2 =\xi(s_p)+s_p -2\overset{\eqref{distanceA}}{=} \xi(s_{p-1})+s_{p-1} -2= 2j_{s_{p-1}}'.
\end{gather*}
It is now immediate to check that \eqref{notturningpoints} holds in this case. The case when $\epsilon_p=1$ is similar and we omit details.\\\\ 
 We now prove the reverse inclusion. If $\epsilon_p=0$ we show that 
  $$\{[i_{s_p}',j] \in \mathbb I_N(\bos_{\xi_\bos}): j\in \{j_{s_{p-1}}', j_{s_{p}}',j_{s_{p+1}}',j_{s_{p+2}}'\}\}\subset \mathbb I_N(\xi_\bos).$$
  This follows from an inspection, noting that $[i_{s_p}',j_{s_{p+2}}'] = 
      [i_{s_{p+1}}^\xi, j_{s_{p+1}}^\xi]$, and 
  \begin{gather*}
  [i_{s_p}',j_{s_{p-1}}'] = \begin{cases} [i_{s_{p}-1}^\xi, j_{s_p-1}^\xi], & {\rm if} \ s_p-1\ne s_{p-1},\\
  [i_{s_{p+1}}^\xi, j_{s_{p+1}}^\xi], & {\rm if} \ s_{p-1}=s_p-1,
  \end{cases}\\ 
  [i_{s_p}',j_{s_{p+1}}'] = \begin{cases}
      [i_{s_{p+1}-1}^\xi-1, j_{s_{p+1}-1}^\xi-1], & {\rm if} \ s_{p+1}-1\ne s_p,\\
      [i_{s_{p+1}}^\xi, j_{s_{p+1}}^\xi], & {\rm if} \ s_{p+1}-1=s_p,
  \end{cases}
  \end{gather*}
Similarly, if $\epsilon_p=1$ one shows that  $$\{[i_{s_p}',j] \in \mathbb I_N(\bos_{\xi_\bos}): j\in \{j_{s_{p-2}}', j_{s_{p-1}}',j_{s_{p}}',j_{s_{p+1}}'\}\}\subset \mathbb I_N(\xi_\bos),$$
  which completes the proof of Proposition \ref{red}.

 \end{document}